\documentclass{amsart}

\usepackage{amsfonts,amssymb,verbatim,amsmath,amsthm,latexsym,textcomp,amscd}
\usepackage{latexsym,amsfonts,amssymb,epsfig,verbatim}
\usepackage{amsmath,amsthm,amssymb,latexsym,graphics,textcomp}
\usepackage{graphicx}
\usepackage{color}
\usepackage{url}

\input xy
\xyoption{all}

\begin{document}

\newtheorem{theorem}{Theorem}[section]
\newtheorem{prop}[theorem]{Proposition}
\newtheorem{lemma}[theorem]{Lemma}
\newtheorem{cor}[theorem]{Corollary}
\newtheorem{definition}[theorem]{Definition}
\newtheorem{conj}[theorem]{Conjecture}
\newtheorem{rmk}[theorem]{Remark}
\newtheorem{claim}[theorem]{Claim}
\newtheorem{defth}[theorem]{Definition-Theorem}

\newcommand{\boundary}{\partial}
\newcommand{\C}{{\mathbb C}}
\newcommand{\integers}{{\mathbb Z}}
\newcommand{\natls}{{\mathbb N}}
\newcommand{\ratls}{{\mathbb Q}}
\newcommand{\bbR}{{\mathbb R}}
\newcommand{\proj}{{\mathbb P}}
\newcommand{\lhp}{{\mathbb L}}
\newcommand{\tube}{{\mathbb T}}
\newcommand{\cusp}{{\mathbb P}}
\newcommand\AAA{{\mathcal A}}
\newcommand\BB{{\mathcal B}}
\newcommand\CC{{\mathcal C}}
\newcommand\DD{{\mathcal D}}
\newcommand\EE{{\mathcal E}}
\newcommand\FF{{\mathcal F}}
\newcommand\GG{{\mathcal G}}
\newcommand\HH{{\mathcal H}}
\newcommand\II{{\mathcal I}}
\newcommand\JJ{{\mathcal J}}
\newcommand\KK{{\mathcal K}}
\newcommand\LL{{\mathcal L}}
\newcommand\MM{{\mathcal M}}
\newcommand\NN{{\mathcal N}}
\newcommand\OO{{\mathcal O}}
\newcommand\PP{{\mathcal P}}
\newcommand\QQ{{\mathcal Q}}
\newcommand\RR{{\mathcal R}}
\newcommand\SSS{{\mathcal S}}
\newcommand\TT{{\mathcal T}}
\newcommand\UU{{\mathcal U}}
\newcommand\VV{{\mathcal V}}
\newcommand\WW{{\mathcal W}}
\newcommand\XX{{\mathcal X}}
\newcommand\YY{{\mathcal Y}}
\newcommand\ZZ{{\mathcal Z}}
\newcommand\CH{{\CC\HH}}
\newcommand\PEY{{\PP\EE\YY}}
\newcommand\MF{{\MM\FF}}
\newcommand\RCT{{{\mathcal R}_{CT}}}
\newcommand\PMF{{\PP\kern-2pt\MM\FF}}
\newcommand\FL{{\FF\LL}}
\newcommand\PML{{\PP\kern-2pt\MM\LL}}
\newcommand\GL{{\GG\LL}}
\newcommand\Pol{{\mathcal P}}
\newcommand\half{{\textstyle{\frac12}}}
\newcommand\Half{{\frac12}}
\newcommand\Mod{\operatorname{Mod}}
\newcommand\Area{\operatorname{Area}}
\newcommand\ep{\epsilon}
\newcommand\hhat{\widehat}
\newcommand\Proj{{\mathbf P}}
\newcommand\U{{\mathbf U}}
 \newcommand\Hyp{{\mathbf H}}
\newcommand\D{{\mathbf D}}
\newcommand\Z{{\mathbb Z}}
\newcommand\R{{\mathbb R}}
\newcommand\Q{{\mathbb Q}}
\newcommand\E{{\mathbb E}}
\newcommand\til{\widetilde}
\newcommand\length{\operatorname{length}}
\newcommand\tr{\operatorname{tr}}
\newcommand\gesim{\succ}
\newcommand\lesim{\prec}
\newcommand\simle{\lesim}
\newcommand\simge{\gesim}
\newcommand{\simmult}{\asymp}
\newcommand{\simadd}{\mathrel{\overset{\text{\tiny $+$}}{\sim}}}
\newcommand{\ssm}{\setminus}
\newcommand{\diam}{\operatorname{diam}}
\newcommand{\pair}[1]{\langle #1\rangle}
\newcommand{\T}{{\mathbf T}}
\newcommand{\inj}{\operatorname{inj}}
\newcommand{\pleat}{\operatorname{\mathbf{pleat}}}
\newcommand{\short}{\operatorname{\mathbf{short}}}
\newcommand{\vertices}{\operatorname{vert}}
\newcommand{\collar}{\operatorname{\mathbf{collar}}}
\newcommand{\bcollar}{\operatorname{\overline{\mathbf{collar}}}}
\newcommand{\I}{{\mathbf I}}
\newcommand{\tprec}{\prec_t}
\newcommand{\fprec}{\prec_f}
\newcommand{\bprec}{\prec_b}
\newcommand{\pprec}{\prec_p}
\newcommand{\ppreceq}{\preceq_p}
\newcommand{\sprec}{\prec_s}
\newcommand{\cpreceq}{\preceq_c}
\newcommand{\cprec}{\prec_c}
\newcommand{\topprec}{\prec_{\rm top}}
\newcommand{\Topprec}{\prec_{\rm TOP}}
\newcommand{\fsub}{\mathrel{\scriptstyle\searrow}}
\newcommand{\bsub}{\mathrel{\scriptstyle\swarrow}}
\newcommand{\fsubd}{\mathrel{{\scriptstyle\searrow}\kern-1ex^d\kern0.5ex}}
\newcommand{\bsubd}{\mathrel{{\scriptstyle\swarrow}\kern-1.6ex^d\kern0.8ex}}
\newcommand{\fsubeq}{\mathrel{\raise-.7ex\hbox{$\overset{\searrow}{=}$}}}
\newcommand{\bsubeq}{\mathrel{\raise-.7ex\hbox{$\overset{\swarrow}{=}$}}}
\newcommand{\tw}{\operatorname{tw}}
\newcommand{\base}{\operatorname{base}}
\newcommand{\trans}{\operatorname{trans}}
\newcommand{\rest}{|_}
\newcommand{\bbar}{\overline}
\newcommand{\UML}{\operatorname{\UU\MM\LL}}
\newcommand{\EL}{\mathcal{EL}}
\newcommand{\tsum}{\sideset{}{'}\sum}
\newcommand{\tsh}[1]{\left\{\kern-.9ex\left\{#1\right\}\kern-.9ex\right\}}
\newcommand{\Tsh}[2]{\tsh{#2}_{#1}}
\newcommand{\qeq}{\mathrel{\approx}}
\newcommand{\Qeq}[1]{\mathrel{\approx_{#1}}}
\newcommand{\qle}{\lesssim}
\newcommand{\Qle}[1]{\mathrel{\lesssim_{#1}}}
\newcommand{\simp}{\operatorname{simp}}
\newcommand{\vsucc}{\operatorname{succ}}
\newcommand{\vpred}{\operatorname{pred}}
\newcommand\fhalf[1]{\overrightarrow {#1}}
\newcommand\bhalf[1]{\overleftarrow {#1}}
\newcommand\sleft{_{\text{left}}}
\newcommand\sright{_{\text{right}}}
\newcommand\sbtop{_{\text{top}}}
\newcommand\sbot{_{\text{bot}}}
\newcommand\sll{_{\mathbf l}}
\newcommand\srr{_{\mathbf r}}
\newcommand\geod{\operatorname{\mathbf g}}
\newcommand\vol{\operatorname{\mathrm vol}}
\newcommand\mtorus[1]{\boundary U(#1)}
\newcommand\A{\mathbf A}
\newcommand\Aleft[1]{\A\sleft(#1)}
\newcommand\Aright[1]{\A\sright(#1)}
\newcommand\Atop[1]{\A\sbtop(#1)}
\newcommand\Abot[1]{\A\sbot(#1)}
\newcommand\boundvert{{\boundary_{||}}}
\newcommand\storus[1]{U(#1)}
\newcommand\Momega{\omega_M}
\newcommand\nomega{\omega_\nu}
\newcommand\twist{\operatorname{tw}}
\newcommand\modl{M_\nu}
\newcommand\MT{{\mathbb T}}
\newcommand\Teich{{\mathcal T}}
\renewcommand{\Re}{\operatorname{Re}}
\renewcommand{\Im}{\operatorname{Im}}

\title[Fourier transform on harmonic manifolds]{The Fourier transform on harmonic manifolds of purely exponential volume growth}

\author{Kingshook Biswas, Gerhard Knieper \and Norbert Peyerimhoff}
\address{Indian Statistical Institute, Kolkata, India. Email: kingshook@isical.ac.in}
\address{Ruhr University Bochum, Germany. Email: gerhard.knieper@rub.de}
\address{Durham University, United Kingdom. Email: norbert.peyerimhoff@durham.ac.uk}

\begin{abstract} Let $X$ be a complete, simply connected harmonic manifold of purely exponential volume growth. This class contains all non-flat harmonic
manifolds of non-positive curvature and, in particular all known examples of harmonic manifolds except for the flat spaces.

Denote by $h > 0$ the mean curvature of horospheres in $X$, and set $\rho = h/2$.
Fixing a basepoint $o \in X$, for $\xi \in \partial X$, denote by $B_{\xi}$ the
Busemann function at $\xi$ such that $B_{\xi}(o) = 0$. then for $\lambda \in \C$
the function $e^{(i\lambda - \rho)B_{\xi}}$ is an eigenfunction of the Laplace-Beltrami operator with
eigenvalue $-(\lambda^2 + \rho^2)$.

For a function $f$ on $X$, we define the Fourier transform of $f$ by
$$\tilde{f}(\lambda, \xi) := \int_X f(x) e^{(-i\lambda - \rho)B_{\xi}(x)} dvol(x)$$ for all
$\lambda \in \C, \xi \in \partial X$ for which the integral converges. We prove a Fourier inversion
formula $$f(x) = C_0 \int_{0}^{\infty} \int_{\partial X} \tilde{f}(\lambda, \xi) e^{(i\lambda - \rho)B_{\xi}(x)} d\lambda_o(\xi) |c(\lambda)|^{-2} d\lambda$$
for $f \in C^{\infty}_c(X)$, where $c$ is a certain function on $\mathbb{R} - \{0\}$, $\lambda_o$ is the
visibility measure on $\partial X$ with respect to the basepoint $o \in X$ and $C_0 > 0$ is a constant. We also prove a Plancherel theorem, and
a version of the Kunze-Stein phenomenon.
\end{abstract}

\bigskip

\maketitle

\tableofcontents

\section{Introduction}

\medskip

Throughout this article, we assume that all manifolds are complete.
A {\it harmonic manifold} is a Riemannian manifold $X$ such that for any point $x \in X$, there exists a
non-constant harmonic function on a punctured neighbourhood of $x$ which is radial around $x$, i.e. only depends
on the geodesic distance from $x$. Copson and Ruse showed that this is equivalent to requiring that sufficiently small
geodesic spheres centered at $x$ have constant mean curvature, and moreover such manifolds are Einstein manifolds \cite{copsonruse}.
Hence they have constant curvature in dimensions 2 and 3.
The Euclidean spaces and rank one symmetric
spaces are examples of harmonic manifolds. The Lichnerowicz conjecture asserts that conversely any harmonic manifold is
either flat or locally symmetric of rank one. The conjecture was proved for harmonic manifolds of dimension $4$ by A. G. Walker \cite{walker1}.
In 1990 Z. I. Szabo proved the conjecture for compact simply connected harmonic manifolds \cite{szabo}.
In 1995 G. Besson, G. Courtois and S. Gallot proved the conjecture for manifolds of negative curvature admitting a
compact quotient \cite{bcg1}, using rigidity results from hyperbolic dynamics including the work of
Y. Benoist, P. Foulon and F. Labourie \cite{bfl} and that of P. Foulon and F. Labourie \cite{foulonlabourie}. In 2005
Y. Nikolayevsky proved the conjecture for harmonic manifolds of dimension 5, showing that these must in fact
have constant curvature \cite{nikolayevsky}. Another fundamental result states that harmonic manifolds of subexponential volume growth are flat \cite{RS02}.

\medskip

In 1992 however E. Damek and F. Ricci had already provided in the non-compact case a family of counterexamples to the Lichnerowicz conjecture, which
have come to be known as {\it harmonic NA groups}, or {\it Damek-Ricci spaces} \cite{damekricci1}. These are solvable Lie groups $X = NA$ with a suitable
left-invariant Riemannian metric, given by
the semi-direct product of a nilpotent Lie group $N$ of {\it Heisenberg type} (see \cite{kaplan}) with $A = \R^+$ acting
on $N$ by anisotropic dilations. While the non-compact rank one symmetric spaces $G/K$ may be identified with harmonic $NA$ groups
(apart from the real hyperbolic spaces), there are examples of harmonic $NA$ groups which are not symmetric.
In 2006, J. Heber proved that the only complete simply connected homogeneous harmonic manifolds are the Euclidean spaces, rank one symmetric spaces,
and harmonic $NA$ groups \cite{heber}.

\medskip

Though the harmonic $NA$ groups are not symmetric in general, there is still a well developed theory of
harmonic analysis on these spaces which parallels that of the symmetric spaces $G/K$. For a non-compact symmetric space
$X = G/K$, an important role in the analysis on these spaces is played by the well-known {\it Helgason Fourier transform}
\cite{helgason1}. For harmonic $NA$ groups, F. Astengo, R. Camporesi and B. Di Blasio have defined a Fourier transform
\cite{astcampblas}, which reduces to the Helgason Fourier transform when the space is symmetric. In both cases a Fourier
inversion formula and a Plancherel theorem hold.

\medskip

The aim of the present article is to generalize these results to
a large class of non-compact harmonic manifolds.
Our analysis will be concerned with harmonic manifolds of purely exponential volume growth which include all non-flat harmonic manifolds
of non-positive sectional curvature or, more generally, all non-flat harmonic manifolds without focal points (see \cite[Theorem 6.5]{knieper1}). In particular this class 
includes all known examples of non-flat and non-compact harmonic manifolds. By {\it purely exponential volume growth}, we mean that there are constants $C > 1$, $h > 0$ such that for all $R > 1$ the volume of metric balls $B(x, R)$ of radius $R$ and center $x \in X$ is given by
\begin{equation} \label{eq:purexpgrowth}
\frac{1}{C} e^{hR} \leq vol(B(x, R)) \leq C e^{hR}.
\end{equation}

Let $X$ be a simply connected harmonic manifold of purely exponential volume growth with a fixed basepoint $o \in X$.
It was shown in \cite{knieper1} that for harmonic manifolds the condition of purely exponential volume growth is equivalent to Gromov hyperbolicity. Moreover, it follows from the work in \cite{knieperpeyerimhoff2} that the Gromov boundary agrees with the visibility boundary $\partial X$ introduced in \cite{eboneill73}. The set $X \cup \partial X$ equipped with the cone topology defines a topological space homeomorphic to a closed unit ball in $\R^n$, where $n = \dim X$. For a given $\xi \in \partial X$ and any geodesic ray $\gamma: [0,\infty) \to X$ representing $\xi$
(see section \ref{sec:basics} for a precise definition)
the Busemann function $B_\xi$ with $B_\xi(o) = 0$ is given by
$$ B_\xi(y) = \lim_{t \to \infty} (d(y,\gamma(t))-d(o, \gamma(t))). $$
The level sets of $B_\xi$ are called {\it horospheres} in $X$. The manifold $X$, being harmonic,
is also {\it asymptotically harmonic}, i.e. the mean curvature of all horospheres is equal to a constant $h \geq 0$.
If $X$ has purely exponential volume growth then $h$ is positive and agrees with the constant $h$ appearing in
\eqref{eq:purexpgrowth}. An easy computation shows that for $\rho = h/2$ and any $\lambda \in \C$ and $\xi \in \partial X$,
the function $f = e^{(i\lambda - \rho)B_{\xi}}$ is an eigenfunction of the Laplace-Beltrami operator $\Delta$ on $X$ with
eigenvalue $-(\lambda^2 + \rho^2)$.

\medskip

The Fourier transform of a function $f \in C^{\infty}_c(X)$ is then defined to be the function on $\C \times \partial X$ given by
$$
\tilde{f}(\lambda, \xi) = \int_X f(x) e^{(-i\lambda - \rho)B_{\xi}(x)} dvol(x).
$$
When $X$ is a non-compact rank one symmetric space, this reduces to the Helgason Fourier transform.

\medskip

The normalized canonical measure of the unit tangent sphere $T^1_o X$ induced by the Riemannian metric
is denoted by $\theta_o$. The unit tangent sphere $T^1_o X$ is identified with the boundary $\partial X$ via the homeomorphism $pr_o : v \in T^1_o X \mapsto \xi = \gamma_v(\infty) \in \partial X$, where $\gamma_v$ is the unique geodesic ray with $\gamma'_v(0) = v$. Pushing forward the measure $\theta_o$ on $T^1_o X$ by the map $pr_o$ gives a measure on $\partial X$ called the {\it visibility measure}, which we denote by $\lambda_o$. We have the following Fourier inversion formula:

\medskip

\begin{theorem} \label{mainthm1} Let $(X,g)$ be a simply connected, harmonic manifold of purely exponential volume growth. Then there is a constant $C_0 > 0$ and a function $c$ on $\C - \{0\}$ such that
for any $f \in C^{\infty}_c(X)$, we have
$$
f(x) = C_0 \int_{0}^{\infty} \int_{\partial X} \tilde{f}(\lambda, \xi) e^{(i\lambda - \rho)B_{\xi}(x)} d\lambda_o(\xi) |c(\lambda)|^{-2} d\lambda
$$
for all $x \in X$.
\end{theorem}

\medskip

We also obtain a Plancherel formula:

\begin{theorem} \label{mainthm2} Let $(X,g)$ be a simply connected, harmonic manifold of purely exponential volume growth. For any $f, g \in C^{\infty}_c(X)$, we have
$$
\int_X f(x) \overline{g(x)} dvol(x) = C_0 \int_{0}^{\infty} \int_{\partial X} \tilde{f}(\lambda, \xi) \overline{\tilde{g}(\lambda, \xi)} d\lambda_o(\xi)
|c(\lambda)|^{-2} d\lambda.
$$
The Fourier transform extends to an isometry of $L^2(X, dvol)$ into $L^2((0, \infty) \times \partial X, C_0 d\lambda_o(\xi) |c(\lambda)|^{-2} d\lambda)$.
\end{theorem}

\medskip

The function $c$ in the previous two theorems is holomorphic on $\Im \lambda < 0$ and has the following integral representation:

\medskip

\begin{theorem} \label{mainthm_c}
Let $(X,g)$ be a simply connected harmonic manifold of purely exponential volume growth
and $c$ be the $c$-function of the radial hypergroup of $X$. Let $\Im \lambda < 0$. Then we have
$$
c(\lambda) = \int_{\partial X} e^{-2(i\lambda - \rho)(\xi|\eta)_x} d\lambda_x(\eta).
$$
for any $x \in X, \xi \in \partial X$, where $(\xi|\eta)_x$ is the Gromov product on $X$ given in Definition \ref{def:gromov_product} below.
\end{theorem}

\medskip

We define a notion of convolution with radial functions and prove the following version of the
Kunze-Stein phenomenon:

\begin{theorem} \label{mainthm3} Let $(X,g)$ be a simply connected harmonic manifold of purely exponential volume growth. Let $x \in X$ and let $1 \leq p < 2$. Let $g \in C^{\infty}_c(X)$ be
radial around the point $x$. Then for any $f \in C^{\infty}_c(X)$ the inequality
$$
||f * g||_2 \leq C_p ||g||_p ||f||_2
$$
holds for some constant $C_p > 0$.
It follows that for any $g \in L^p(X)$ radial around $x$, the map $f \in C^{\infty}_c(X) \mapsto f * g$
extends to a bounded linear operator on $L^2(X)$ with operator norm at most $C_p ||g||_p$.
\end{theorem}

\medskip

The article is organized as follows. In section 2 we recall basic facts about harmonic manifolds
which we require. In section 3 we compute the action of the Laplacian $\Delta$ on spaces of functions constant on
geodesic spheres and horospheres respectively. In section 4 we carry out the harmonic analysis of radial functions,
i.e. functions constant on geodesic spheres centered around a given point. Unlike the well-known {\it Jacobi analysis}
\cite{koornwinder}
which applies to radial functions on rank one symmetric spaces and harmonic $NA$ groups, our analysis here is based
on {\it hypergroups} \cite{bloomheyer}. We define a spherical Fourier transform for radial functions, and obtain an
inversion formula and Plancherel theorem for this transform. In section 5 we prove the inversion formula
and Plancherel formula for the Fourier transform. The main point of the proof is an identity expressing
radial eigenfunctions in terms of an integral over the boundary $\partial X$. The integral formula for the function $c$ (Theorem \ref{mainthm_c}) is proved in section 6. In section 7 we define
an operation of convolution with radial functions, and show that the $L^1$ radial functions form a commutative
Banach algebra under convolution. Finally in section 8 we prove a version of the Kunze-Stein phenomenon.

\medskip

{\bf Acknowledgements.} The first author would like to thank Swagato K. Ray and Rudra P. Sarkar for generously
sharing their time and knowledge over the course of numerous educative and enjoyable discussions. The other two authors like to thank the MFO for hospitality during their stay in the "Research in Pairs" program in 2019 and the SFB/TR191 "Symplectic structures in geometry, algebra and dynamics". This article generalizes an earlier version by the first author in the case of negatively curved harmonic manifolds.

\medskip

\section{Basics about harmonic manifolds}\label{sec:basics}

Throughout this article, we assume that all manifolds are complete.
We start by presenting some fundamental facts about  non-compact simply connected harmonic manifolds. References for this class of manifolds include \cite{rusewalkerwillmore}, \cite{szabo}, \cite{willmore}, \cite{knieperpeyerimhoff1} and \cite{knieper3}.
Such manifolds do not have conjugate points and, for every $x \in X$, the exponential map $\exp_x: T_x X \to X$ is a diffeomorphism. (See e.g \cite{knieper2} on basic geometric  and dynamical properties of spaces without conjugate points.) The absence of conjugate points in $X$ allows to define Busemann functions associated to geodesic rays $\gamma_v: [0,\infty) \to X$ with $\gamma_v'(0) = v$. These functions are of central importance in our paper and are given by
$$ b_v(y) = \lim_{t \to \infty} (d(y,\gamma_v(t))-t). $$
The level sets of these functions are called {\it horospheres} and can be viewed as spheres with center at infinity.

For any $v \in T^1_x X$ and $r > 0$, let $A(v, r)$ denote the Jacobian of the map $v \mapsto \exp_x(rv)$. The definition of harmonicity given in the Introduction is equivalent to the fact that this Jacobian does not depend on $v$, i.e. there is a
function $A$ on $(0, \infty)$ such that $A(v, r) = A(r)$ for all $v \in T^1 X$. See \cite[p. 224]{willmore} for the equivalence of this property with the property given in the Introduction. The function $A$ is called the
{\it density function} of the harmonic manifold.

For $x \in X$, let $d_x$ denote the distance function from the point $x$, i.e. $d_x(y) = d(x, y)$.
A function $f$ on $X$ is said to be {\it radial} around a point $x$ of $X$ if $f$ is constant on geodesic
spheres centered at $x$. For each $x \in X$, we can define a radialization operator $M_x$, defined for a
continuous function $f$ on $X$ by
$$
(M_x f)(z) = \int_{S(x, r)} f(y) d\sigma(y)
$$
where $S(x, r)$ denotes the geodesic sphere around $x$ of radius $r = d(x, z)$, and $\sigma$ denotes surface
area measure on this sphere (induced from the metric on $X$), normalized to have mass one.
The operator $M_x$ maps continuous functions to functions radial around $x$, and is formally self-adjoint,
meaning
$$
\int_X (M_x f)(z) h(z) dvol(z) = \int_X f(z) (M_x h)(z) dvol(z).
$$
for all continuous functions $f, h$ with compact support. Introducing polar coordinates around $x$ this follows
easily from
$$
\int_X (M_x f)(z) h(z) dvol(z) = \int_0^\infty \int_{T^1_xX} f(\gamma_v(r)) d\theta_x(v) \int_{T^1_xX} h(\gamma_w(r)) d\theta_x(w)  A(r) dr,
$$
where $\theta_x$ is the normalized canonical measure on the unit tangent space $T^1_xX$ induced by the Riemannian metric and $\gamma_v: \R \to X$ is the geodesic satisfying $\gamma_v'(0)=v$.

Using these concepts, we have the following equivalent conditions for harmonicity:
\begin{itemize}
\item[(1)] For any $x \in X$, $\Delta d_x$ is radial around $x$.
\item[(2)] The Laplacian $\Delta = {\rm{div}} \circ \nabla$ commutes with all the radialization operators $M_x$, i.e. $M_x \Delta u = \Delta M_x u$ for all smooth functions $u$ on $X$ and all $x \in X$.
\item[(3)] For any smooth function $u$ radial around any $x \in X$ the function $\Delta u$ is radial around $x$, as well.
\end{itemize}

Let us now discuss basic properties of the density function $A(r)$ of a harmonic manifold. $A(r)$ is increasing in $r$, and the quantity $A'(r)/A(r) \geq 0$ equals the mean curvature of geodesic spheres $S(x, r)$ of radius $r$, which decreases monotonically as $r \to \infty$ (see \cite[Corollary 2.1, Proposition 2.2]{RS03} and \cite[Section 1.2]{knieper2}). Furthermore, the mean curvature $(A'/A)(r)$ of the geodesic sphere $S(x, r)$ at a point $z \in S(x, r)$ equals $\Delta d_x (z)$, hence
we have
$$
\Delta d_x = \frac{A'}{A} \circ d_x.
$$
The limit $\lim_{r \to \infty} A'(r)/A(r)$ is equal to the mean curvature $h \ge 0$ of horospheres. Therefore, all harmonic manifolds are in particular \emph{asymptotically harmonic}, meaning they are manifolds without conjugate points such that all horospheres have the same constant mean curvature.

Using the density function $A(r)$, harmonic manifolds are of purely exponential volume growth if and only if there
exist constants $C > 1$, $h > 0$ such that we have for all $R > 1$
$$ \frac{1}{C} e^{h R} \le A(R) \le C e^{h R}. $$
In this particular case it turns out that the constant $h > 0$ agrees with the mean curvature of the horospheres.

Let us finish this section by discussing specific properties of non-compact simply connected \emph{harmonic manifolds} $(X,g)$ of \emph{purely exponential volume growth} as defined in \eqref{eq:purexpgrowth}.  In this setting, purely exponential volume growth, Anosov geodesic flow and Gromov hyperbolicity are equivalent properties (see \cite{knieper1}). A geodesic metric space $(X,d)$ is called \emph{Gromov hyperbolic} if there exists a $\delta > 0$ such that geodesic triangles are $\delta$-thin,
that is each side is contained in the $\delta$-tubes of the other two sides. 

\medskip

Next we introduce a boundary structure for $(X,g)$ and define a natural topology. The boundary structure is given by  equivalence classes of geodesic rays in $X$, where two rays $\gamma_1, \gamma_2$ are equivalent if $\{d(\gamma_1(t),\gamma_2(t)): t \ge 0 \}$ is bounded. We denote this boundary by $\partial X$ and the equivalence class associated to a geodesic ray $\gamma$ by $\gamma(\infty) \in \partial X$.
Let $\bar X = X \cup \partial X$. For each $x \in X$, we introduce the following bijective map $pr_x: B_1(x) \to \bar X$, where $B_1(x) \subset T_xX$ is the closed ball of radius $1$:
$$ pr_x(v) = \begin{cases} \gamma_v(\infty) & \text{if $\Vert v \Vert = 1$,} \\ \exp_x(\frac{1}{1-\Vert v \Vert} v
& \text{if $\Vert v \Vert < 1$.} \end{cases} $$
Then the topology on $\bar X$ is defined such that $pr_x$ is a homeomorphism. This definition does not depend on the choice of $x$ and is called the {\it cone topology}. We proved in \cite[Theorem 4.5]{knieperpeyerimhoff2} that this topology agrees with the Gromov topology on $\bar X$. 

Since the horospheres are the footpoint projections of the stable manifolds of the geodesic flow, we have the following convergence property of asymptotic geodesic starting from the same horosphere in the case of Anosov geodesic flow: 
given $\xi = \gamma_v(\infty) \in \partial X$ and $x, y \in X$ such that $b_v(x) = b_v(y) = 0$, and geodesics $\gamma_1, \gamma_2 : [0, +\infty) \to X$ such that $\gamma_1(0) = x, \gamma_2(0) = y$
and $\gamma_1(\infty) = \gamma_2(\infty) = \xi$, we have that $d(\gamma_1(t), \gamma_2(t)) \to 0$ as $t \to \infty$.
Using this fact we define Busemann functions alternatively with respect to boundary points as follows:

\begin{lemma} \label{busemann_intro}
  Let $(X,g)$ be a simply connected harmonic manifold of purely exponential volume growth and $x \in X$ and $\xi \in \partial X$. Then the \emph{Busemann function} $B_{\xi,x}: X \to \R$ is defined by
  $$ B_{\xi,x}(y) = \lim_{t \to \infty} (d(y,\gamma(t))-d(x,\gamma(t))) $$
where $\gamma: [0,\infty) \to X$ is a geodesic ray with $\gamma(\infty) = \xi$. This definition does not depend
on the choice of $\gamma$.
\end{lemma}

\medskip

\noindent{\bf Proof:} Let $\gamma_0: [0,\infty) \to X$ be the geodesic ray with $\gamma_0(0) = x$
and $\gamma_0(\infty) = \xi$. Let $v = \gamma_0'(0)$. Then there exists $t_0 \in \R$ such that we have
$$ d(\gamma_0(t+t_0),\gamma(t)) \to 0 \quad \text{for $t \to \infty$,} $$
and we have
\begin{multline*}
   d(y,\gamma(t))-d(x,\gamma(t)) = d(y,\gamma_0(t+t_0)) + ( d(y,\gamma(t)) - d(y,\gamma_0(t+t_0)) )\\ - d(x,\gamma_0(t+t_0)) - (  d(x,\gamma(t)) - d(x,\gamma_0(t+t_0)) ). 
\end{multline*}
Since 
\begin{equation} \label{eq:triangl-useful} 
| d(z,\gamma(t)) - d(z,\gamma_0(t+t_0)) | \le d(\gamma(t),\gamma_0(t+t_0)) \to 0 \quad \text{for $t \to \infty$,} 
\end{equation}
we obtain
\begin{multline*} 
\lim_{t \to \infty} (d(y,\gamma(t))-d(x,\gamma(t))) = \lim_{t \to \infty} (d(y,\gamma_0(t+t_0)) - d(x,\gamma_0(t+t_0))) 
= \\ \lim_{t \to \infty} (d(y,\gamma_0(t+t_0)) - (t+t_0) = b_v(q). 
\end{multline*}
This shows the independence of the limit of the choice of geodesic ray.
$\diamond$

\medskip

The level sets of $B_{\xi,x}$ are called {\it horospheres} centered at $\xi$ and their mean curvatures agree with  $\Delta B_{\xi, x}$ for all $\xi \in \partial X, x \in X$.
Since they have the same constant mean curvature $h \ge 0$, we have
$$ \Delta B_{\xi, x} = h. $$
In the case of purely exponential volume growth the constant $h$ is positive.
The \emph{Busemann cocycle} $B: \partial X \times X \times X \to \R$ is defined by
$$ B(x,y,\xi) := B_{\xi,y}(x), $$
and it is easy to see that it satisfies the following cocycle property:
$$ B(x,z,\xi) = B(x,y,\xi) + B(y,z,\xi). $$

\medskip

Since $(X,g)$ is a Gromov hyperbolic space by \cite{knieper1}, it is equipped with the Gromov product defined as follows (see \cite{BS2007}):

\begin{lemma} \label{def:gromov_product}
Let $(X,g)$ be a simply connected harmonic manifold of purely exponential volume growth.
For given $x \in X$ and any $y,z \in X$
we define
$$ (y|z)_x = \frac{1}{2}\left( d(x,y) + d(x,z) - d(,y,z) \right) \ge 0, $$
and for $\xi,\eta \in \partial X$,
$$
(\xi|\eta)_x = \lim_{s,t \to \infty} (\gamma_1(s)|\gamma_2(t))_x,
$$
where $\gamma_1,\gamma_2: [0,\infty) \to X$ are geodesic rays with $\gamma_1(\infty) = \xi$
and $\gamma_2(\infty) = \eta$. This definition does not depend on the choice of $\gamma_1,\gamma_2$.
The map $(\cdot|\cdot)_x: \bar X \times \bar X \to [0, +\infty]$ is called the \emph{Gromov product}.
\end{lemma}

\medskip

\noindent{\bf Proof:} We first assume $\xi \neq \eta$. Since $X$ is Gromov hyperbolic, there exists a geodesic
$\gamma: \R \to X$ with $\gamma(-\infty) = \xi$ and $\gamma(\infty) = \eta$
 (see, e.g., \cite[Lemma 11.83]{DK18}). Using the Anosov property, we conclude that there exist $s_0,t_0 \in \R$ such that
$$ d(\gamma_1(s),\gamma(-s+s_0)) \to 0 \quad \text{as $s \to \infty$} $$
and 
$$ d(\gamma_2(t),\gamma(t+t_0)) \to 0 \quad \text{as $t \to \infty$}. $$
Using these limits and similar arguments as in the proof of Lemma \ref{busemann_intro} (in particular \eqref{eq:triangl-useful}), we derive
\begin{multline*}
\lim_{s,t \to \infty} (\gamma_1(s)|\gamma_2(t))_x = 
\lim_{s,t \to \infty} \frac{1}{2} ( d(\gamma_1(s),x) + d(\gamma_2(t),x) - d(\gamma_1(s),\gamma_2(t)) ) = \\
\lim_{s,t \to \infty} \frac{1}{2} ( d(\gamma(-s+s_0),x) + d(\gamma(t+t_0),x) - d(\gamma(-s+s_0),\gamma(t+t_0)) ) = \\
\frac{1}{2} ( \lim_{s \to \infty} ( d(\gamma(-s+s_0),x) - (s-s_0) ) ) + \frac{1}{2} ( \lim_{t \to \infty} ( d(\gamma(t+t_0),x) - (t+t_0) ) ) = \\
\frac{1}{2} ( B_{\xi,\gamma(0)}(x) + B_{\eta,\gamma(0)}(x) ).
\end{multline*}
Next we assume $\xi = \eta$. Let $\gamma, \gamma_1, \gamma_2: [0,\infty) \to X$ be geodesic rays with
$\gamma(0) = x$ and $\gamma(\infty) = \gamma_1(\infty) = \gamma_2(\infty) = \xi$. Again, we can find $t_1, t_2 \in \R$ such that for $i \in \{1,2\}$,
$$ d(\gamma_i(t),\gamma(t+t_i)) \to 0 \quad \text{as $t \to \infty$}. $$
Using these limits again we derive
\begin{multline*}
\lim_{s,t \to \infty} (\gamma_1(s)|\gamma_2(t))_x = 
\lim_{s,t \to \infty} \frac{1}{2} ( d(\gamma_1(s),x) + d(\gamma_2(t),x) - d(\gamma_1(s),\gamma_2(t)) ) = \\
\lim_{s,t \to \infty} \frac{1}{2} ( d(\gamma(s+t_1),x) + d(\gamma(t+t_2),x) - d(\gamma(s+t_1),\gamma(t+t_2)) ) = \\
\lim_{s,t \to \infty} \frac{1}{2} ( s + t_1 + t + t_2 - | s + t_1 - (t + t_2) | ) = \infty.
\end{multline*} 
$\diamond$

\medskip

We have the following relation between Busemann functions and the Gromov product in our setting (it also holds in any CAT(-1) space):

\medskip

\begin{lemma} \label{buseform}
Let $X$ be a noncompact, simply connected harmonic manifold of purely exponential volume growth.
For $x \in X$
and $\eta \in \partial X$, let $\gamma_{x,\eta}: [0,\infty) \to X$ be a geodesic ray with $\gamma_{x,\eta}(0) = x$
and $\gamma_{x,\eta}(\infty) = \eta$. Then we have for all $\xi \in \partial X$:
$$ \lim_{r \to \infty} \left( B_{\xi,x}(\gamma_{x,\eta}(r)) - r \right) = -2 (\xi|\eta)_x. $$
\end{lemma}

\medskip

\noindent{\bf Proof:} Let $\alpha : [0, \infty) \to X$ be a geodesic ray with $\alpha(0) = x$ and $\alpha(\infty) = \xi$.
Then by the previous Lemma, the double limit
$$
\lim_{s, r \to \infty} d(\alpha(s), \gamma_{x, \eta}(r)) - (r + s)
$$
exists and equals $-2(\xi|\eta)_x$. Since the double limit exists, it can be evaluated as an iterated limit, so we have:
$$
-2(\xi|\eta)_x = \lim_{r \to \infty} \left(\lim_{s \to \infty} d(\alpha(s), \gamma_{x, \eta}(r)) - (r + s)\right)
$$
Now for a fixed $r$ we have $\lim_{s \to \infty} (d(\alpha(s), \gamma_{x, \eta}(r)) - (r + s)) = B_{\xi, x}(\gamma_{x,\eta}(r)) - r$, so
substituting this in the previous equation gives the result. $\diamond$
%We have
%\begin{eqnarray*}
%B_{\xi,x}(\gamma_{x,\eta}(r)) - r &=& \lim_{z \to \xi} d(\gamma_{x,\eta}(r),z) - d(x,z) - r\\
%&=& - \lim_{z \to \xi} d(x,z) + d(x,\gamma_{x,\eta}(r)) -  d(\gamma_{x,\eta}(r),z) \\
%&=& - 2 \lim_{z \to \xi} (\gamma_{x,\eta}(r)|z)_x.
%\end{eqnarray*}
%The result follows by taking the limit $r \to \infty$ on both sides.
%$\diamond$

\medskip

Finally, we define the family of visibility measures $\lambda_x$ on harmonic manifolds $(X,g)$ of purely exponential volume growth. For $x \in X$, let $\theta_x$ denote the normalized canonical measure on $T^1_xX$ induced by the Riemannian metric and $\lambda_x$ be the push forward of $\theta_x$ to the boundary $\partial X$ under $p_x$.  The {\it visibility measures} $\lambda_x$ are pairwise absolutely continuous with Radon-Nykodym derivative given by
\begin{equation} \label{eq:RadNyk}
\frac{d\lambda_y}{d\lambda_x}(\xi) = e^{-h B_{\xi,x}(y)}.
\end{equation}
This result was shown in \cite[Theorem 1.4]{knieperpeyerimhoff2} in the more general setting of asymptotically harmonic manifolds of purely exponential volume growth with curvature tensor bounds $\Vert R \Vert_\infty \le R_0$, $\Vert \nabla R \Vert_\infty \le R_0'$ for some $R_0,R_0' > 0$. These curvature tensor bounds are satisfied for harmonic manifolds by \cite[Propositions 6.57 and 6.68]{besse1}.

\medskip

\section{Radial and horospherical parts of the Laplacian}

\medskip

 Let $X$ be a non-compact simply connected harmonic manifold. Let $h \ge 0$ denote the
 mean curvature of horospheres in $X$, let $\rho = \frac{1}{2}h$, and let $A : (0, \infty) \to \R$ denote the density
 function of $X$.

\medskip

\begin{lemma} \label{laplacecomposition} For $f$ a $C^2$ function on $X$ and $u$ a $C^{\infty}$ function on
$\R$, we have
$$
\Delta(u \circ f) = (u'' \circ f) |\nabla f|^2 + (u' \circ f) \Delta f
$$
\end{lemma}

\medskip

\noindent{\bf Proof:}  Let $\gamma$ be a geodesic, then
$(u \circ f \circ \gamma)'(t) = (u' \circ f)(\gamma(t))< \nabla f, \gamma'(t) >$,
so
$$
(u \circ f \circ \gamma)''(t) = (u'' \circ f)(\gamma(t))< \nabla f, \gamma'(t) >^2 + (u' \circ f)(\gamma(t)) < \nabla_{\gamma'} \nabla f, \gamma'(t) >
$$
Now let $\{e_i\}$ be an orthonormal basis of $T_x X$, and let $\gamma_i$ be geodesics
with $\gamma'_i(0) = e_i$. Then
\begin{align*}
\Delta(u \circ f)(x) & = \sum_{i = 1}^n < \nabla_{e_i} \nabla (u \circ f), e_i > \\
                     & = \sum_{i = 1}^n (u \circ f \circ \gamma_i)''(0) \\
                     & = (u'' \circ f)(x) \sum_{i = 1}^n < \nabla f, e_i >^2 + (u' \circ f)(x) \sum_{i = 1}^n < \nabla_{e_i} \nabla f, e_i > \\
                     & = (u'' \circ f)(x) |\nabla f(x)|^2 + (u' \circ f)(x) \Delta f (x) \\
\end{align*}
 $\diamond$

\medskip

Any $C^{\infty}$ function on $X$ radial around $x \in X$ is of the form $f = u \circ d_x$ for some even
$C^{\infty}$ function $u$ on $\R$, where $d_x$ denotes the distance function from the point $x$, while
any $C^{\infty}$ function which is constant on horospheres at $\xi \in \partial X$ is of the form
$f = u \circ B_{\xi, x}$ for some $C^{\infty}$ function $u$ on $\R$. The following proposition says that
the Laplacian $\Delta$ leaves invariant these spaces of functions, and describes the action of the Laplacian
on these spaces:

\medskip

\begin{prop} \label{radialhorospherical} Let $x \in X, \xi \in \partial X$.

\medskip

\noindent (1) For $u$ a $C^{\infty}$ function on $(0, \infty)$,
$$
\Delta (u \circ d_x) = (L_R u) \circ d_x
$$
where $L_R$ is the differential operator on $(0, \infty)$ defined by
$$
L_R = \frac{d^2}{dr^2} + \frac{A'(r)}{A(r)} \frac{d}{dr}
$$

\noindent (2) For $u$ a $C^{\infty}$ function on $\R$,
$$
\Delta (u \circ B_{\xi, x}) = (L_H u) \circ B_{\xi, x}
$$
where $L_H$ is the differential operator on $\R$ defined by
$$
L_H = \frac{d^2}{dt^2} + 2\rho \frac{d}{dt}
$$
\end{prop}

\medskip

\noindent{\bf Proof:} Noting that $|\nabla d_x| = 1, |\nabla B_{\xi, x}| = 1$, and
$\Delta d_x = (A' / A) \circ d_x, \Delta B_{\xi,x} = 2\rho$, the Proposition follows
immediately from the previous Lemma. $\diamond$

\medskip

Accordingly, we call the differential operators $L_R$ and $L_H$
the {\it radial and horospherical parts of the Laplacian} respectively. It follows from the above
proposition that a function $f = u \circ d_x$ radial around $x$ is an eigenfunction of $\Delta$
with eigenvalue $\sigma$ if and only if $u$ is an eigenfunction of $L_R$ with eigenvalue $\sigma$.
Similarly, a function $f = u \circ B_{\xi,x}$ constant on horospheres at $\xi$ is an eigenfunction of
$\Delta$ with eigenvalue $\sigma$ if and only if $u$ is an eigenfunction of $L_H$ with eigenvalue $\sigma$.
In particular, we have the following:

\medskip

\begin{prop} \label{helgasonkernel} Let $\xi \in \partial X, x \in X$. Then for any $\lambda \in \C$, the function
$$
f = e^{(i\lambda - \rho)B_{\xi,x}}
$$
is an eigenfunction of the Laplacian with eigenvalue $-(\lambda^2 + \rho^2)$ satisfying $f(x) = 1$.
\end{prop}

\medskip

\noindent{\bf Proof:} This follows from the fact that
the function $u(t) = e^{(i\lambda - \rho)t}$ on $\R$ is an eigenfunction of $L_H$ with eigenvalue $-(\lambda^2 + \rho^2)$,
and $B_{\xi,x}(x) = 0$ gives $f(x) = 1$. $\diamond$

\medskip

\section{Analysis of radial functions}

\medskip

As we saw in the previous section, finding radial eigenfunctions of the Laplacian
amounts to finding eigenfunctions of its radial part $L_R$. When $X$ is a rank one
symmetric space $G/K$, or more generally a harmonic $NA$ group, then the volume density
function is of the form $A(r) = C \left(\sinh\left(\frac{r}{2}\right)\right)^p \left(\cosh\left(\frac{r}{2}\right)\right)^q$,
for a constant $C > 0$ and integers $p, q \geq 0$, and so the radial part $L_R = \frac{d^2}{dr^2} + (A'/A) \frac{d}{dr}$
falls into the general class of {\it Jacobi operators}
$$
L_{\alpha, \beta} = \frac{d^2}{dr^2} + ((2\alpha + 1)\coth r + (2\beta + 1)\tanh r)\frac{d}{dr}
$$
for which there is a detailed and well known harmonic analysis in terms of eigenfunctions (called {\it Jacobi functions})
\cite{koornwinder}.
For a general harmonic manifold $X$, the explicit form of the density function $A$ is not known, so it is unclear whether
the radial part $L_R$ is a Jacobi operator.
However, there is a harmonic analysis, based on hypergroups
(\cite{chebli1}, \cite{chebli2}, \cite{trimeche1}, \cite{trimeche2}, \cite{trimeche3}, \cite{bloom1}, \cite{xu1}),
for more general second-order differential operators on $(0, \infty)$ of the form
\begin{equation} \label{formofl}
L = \frac{d^2}{dr^2} + \frac{A'(r)}{A(r)} \frac{d}{dr}
\end{equation}
where $A$ is a function on $[0,\infty)$ satisfying certain hypotheses which allow one to endow $[0,\infty)$
with a hypergroup structure, called a {\it Chebli-Trimeche hypergroup}. We first recall some basic facts about Chebli-Trimeche hypergroups,
and then show that the density function of a harmonic manifold satisfies the hypotheses required in order to apply this theory.

\medskip

\subsection{Chebli-Trimeche hypergroups}

\medskip

A hypergroup $(K, *)$ is a locally compact Hausdorff space $K$ such that
the space $M^b(K)$ of finite Borel measures on $K$ is endowed with a product $(\mu, \nu) \mapsto \mu * \nu$
turning it into an algebra with unit, and $K$ is endowed with an involutive homeomorphism $x \in K \mapsto \tilde{x} \in K$,
such that the product and the involution satisfy certain natural properties (see \cite{bloomheyer} Chapter 1 for the
precise definition). A motivating example relevant to the following is the algebra of finite radial measures on a noncompact
rank one symmetric space $G/K$ under convolution; as radial measures can be viewed as measures on $[0, \infty)$, this endows
$[0, \infty)$ with a hypergroup structure (with the involution being the identity).
It turns out that this hypergroup structure on $[0, \infty)$ is a special case
of a general class of hypergroup structures on $[0, \infty)$ called {\it Sturm-Liouville hypergroups} (see \cite{bloomheyer}, section 3.5).
These hypergroups arise from Sturm-Liouville boundary problems on $(0, \infty)$. We will be interested in a particular class
of Sturm-Liouville hypergroups called {\it Chebli-Trimeche hypergroups}. These arise as follows (we refer to \cite{bloomheyer}
for proofs of statements below):

\medskip

A {\it Chebli-Trimeche function} is a continuous function $A$ on $[0, \infty)$ which is $C^{\infty}$ and
positive on $(0, \infty)$ and satisfies the
following conditions:

\medskip

\noindent (H1) $A$ is increasing, and $A(r) \to +\infty$ as $r \to +\infty$.

\medskip

\noindent (H2) $A'/A$ is decreasing, and $\rho = \frac{1}{2} \lim_{r \to \infty} A'(r)/A(r) > 0$.

\medskip

\noindent (H3) For $r > 0$, $A(r) = r^{2\alpha + 1} B(r)$ for some $\alpha > -1/2$ and some even, $C^{\infty}$ function $B$ on $\R$
such that $B(0) = 1$. 

\medskip

Let $L$ be the differential operator on $C^2(0, \infty)$ defined by equation (\ref{formofl}),
where $A$ satisfies conditions (H1)-(H3) above. Define the differential operator $l$ on $C^2((0,\infty)^2)$
by
\begin{align*}
l[u](x,y) & = (L)_x u(x, y) - (L)_y u(x,y) \\
          & = \left(u_{xx}(x,y) + \frac{A'(x)}{A(x)}u_x(x,y)\right) - \left(u_{yy}(x,y) + \frac{A'(y)}{A(y)}u_y(x,y)\right) \\
\end{align*}

For $f \in C^2([0, \infty))$ denote by $u_f$ the solution of the hyperbolic Cauchy problem
\begin{align*}
l[u_f] & = 0, \\
u_f(x, 0) = u_f(0, x) & = f(x), \\
(u_f)_y(x, 0) & = 0, \\
(u_f)_x(0, y) & = 0 \ \hbox{ for } x, y \in [0, \infty) \\
\end{align*}

For $x \in [0, \infty)$, let $\epsilon_x$ denote the Dirac measure of mass one at $x$. Then for all $x, y \in [0, \infty)$, there
exists a probability measure on $[0, \infty)$ denoted by $\epsilon_x * \epsilon_y$ such that
$$
\int_{0}^{\infty} f d(\epsilon_x * \epsilon_y) = u_f(x, y)
$$
for all even, $C^{\infty}$ functions $f$ on $\R$.  We have $\epsilon_x * \epsilon_y = \epsilon_y * \epsilon_x$ for all $x,y$, and
the product
$(\epsilon_x, \epsilon_y) \mapsto \epsilon_x * \epsilon_y$ extends to a product on all finite measures on
$[0, \infty)$ which turns $[0, \infty)$ into a commutative hypergroup $([0,\infty), *)$ (with the involution being the identity),
called the Chebli-Trimeche hypergroup associated to the function $A$.
Any hypergroup has a Haar measure, which in this case
is given by the measure $A(r) dr$ on
$[0, \infty)$.

\medskip

For a commutative hypergroup $K$ with a Haar measure $dk$, a Fourier analysis can be carried out analogous to the Fourier analysis
on locally compact abelian groups. There is a dual space $\hat{K}$ of characters, which are bounded multiplicative functions
on the hypergroup $\chi : K \to \C$ satisfying $\chi(\tilde{x}) = \overline{\chi(x)}$, where multiplicative means that
$$
\int_K \chi d(\epsilon_x * \epsilon_y) = \chi(x) \chi(y)
$$
for all $x, y \in K$. For $f \in L^1(K)$, the Fourier transform of $f$ is the function $\hat{f}$ on $\hat{K}$ defined by
$$
\hat{f}(\chi) = \int_K f \overline{\chi} dk
$$
The Levitan-Plancherel Theorem states that there is a measure $d\chi$ on $\hat{K}$
called the Plancherel measure, such that the mapping $f \mapsto \hat{f}$ extends from
$L^1(K) \cap L^2(K)$ to an isometry from $L^2(K)$ onto $L^2(\hat{K})$. The inverse Fourier transform of a function
$\sigma \in L^1(\hat{K})$ is the function $\check{\sigma}$ on $K$ defined by
$$
\check{\sigma}(k) = \int_{\hat{K}} \sigma(\chi) \chi(k) d\chi
$$
The Fourier inversion theorem then states that if $f \in L^1(K) \cap C(K)$ is such that $\hat{f} \in L^1(\hat{K})$,
then $f = (\hat{f})\check{}$, i.e.
$$
f(x) = \int_{\hat{K}} \hat{f}(\chi) \chi(x) d\chi
$$
for all $x \in K$.

\medskip

For the Chebli-Trimeche hypergroup, it turns out that the multiplicative functions on the hypergroup are
given precisely by eigenfunctions of the operator $L$.
For any $\lambda \in \C$, the equation
\begin{equation} \label{evaleqn}
Lu = -(\lambda^2 + \rho^2)u
\end{equation}
has a unique solution $\phi_{\lambda}$ on $(0, \infty)$ which extends continuously to $0$ and satisfies
$\phi_{\lambda}(0) = 1$ (note that the coefficient $A'/A$ of the operator $L$ is singular at $r = 0$ so existence of a
solution continuous at $0$ is not immediate). The function $\phi_{\lambda}$ extends to a $C^{\infty}$ even function on $\R$.
Since equation (\ref{evaleqn}) reads the same for $\lambda$ and $-\lambda$, by uniqueness we have $\phi_{\lambda} = \phi_{-\lambda}$.

\medskip

The multiplicative functions on $[0, \infty)$ are then exactly the functions $\phi_{\lambda}, \lambda \in \C$. The functions
$\phi_{\lambda}$ are bounded if and only if $| \Im \lambda | \leq \rho$. Furthermore, the involution on the hypergroup being the
identity, the characters of the hypergroup are real-valued, which occurs for $\phi_{\lambda}$ if and only if $\lambda \in \R \cup i\R$.
Thus the dual space of the hypergroup is given by
$$
\hat{K} = \{ \phi_{\lambda} | \lambda \in [0, \infty) \cup [0, i\rho] \}
$$
which we identify with the set $\Sigma = [0, \infty) \cup [0, i\rho] \subset \C$.

\medskip

The hypergroup Fourier transform of a function $f \in L^1([0, \infty), A(r) dr)$ is given by
$$
\hat{f}(\lambda) = \int_{0}^{\infty} f(r) \phi_{\lambda}(r) A(r) dr
$$
for $\lambda \in \Sigma$ (when the hypergroup arises from convolution of
radial measures on a rank one symmetric space $G/K$, then this is the well-known Jacobi transform \cite{koornwinder}).
The Levitan-Plancherel and Fourier inversion theorems for the hypergroup give the existence
of a Plancherel measure $\sigma$ on $\Sigma$ such that the Fourier transform defines an isometry from $L^2([0, \infty), A(r) dr)$ onto
$L^2(\Sigma, \sigma)$, and, for any function $f \in L^1([0, \infty), A(r) dr) \cap C([0, \infty))$ such that $\hat{f} \in L^1(\Sigma, \sigma)$,
we have
$$
f(r) = \int_{\Sigma} \hat{f}(\lambda) \phi_{\lambda}(r) d\sigma(\lambda)
$$
for all $r \in [0, \infty)$.

\medskip

In \cite{bloom1}, it is shown that under certain extra conditions on the function $A$, the support of the Plancherel
measure is $[0, \infty)$ and the Plancherel measure is absolutely continuous with respect to Lebesgue measure $d\lambda$ on $[0,\infty)$,
given by
$$
d\sigma(\lambda) = C_0 |c(\lambda)|^{-2} d\lambda
$$
where $C_0 > 0$ is a constant, and $c$ is a certain complex function on $\C - \{0\}$.
%(when the hypergroup is the one coming from a rank one
%symmetric space, then this function coincides with Harish-Chandra's $c$-function on a closed half-plane, but not on all
%of $\C$)
The required conditions on $A$ are as follows:

\medskip

Making the change of dependent variable $v = A^{1/2} u$, equation (\ref{evaleqn}) becomes
\begin{equation} \label{evaleqn2}
v''(r) = (G(r) - \lambda^2)v(r)
\end{equation}

where the function $G$ is defined by

\begin{equation} \label{gfunction}
G(r) = \frac{1}{4} \left( \frac{A'(r)}{A(r)}\right)^2 + \frac{1}{2}\left(\frac{A'}{A}\right)'(r) - \rho^2
\end{equation}

If the function $G$ tends to $0$ fast enough near infinity, then it is reasonable to expect that equation (\ref{evaleqn2})
above has two linearly independent solutions asymptotic to exponentials $e^{\pm i \lambda r}$ near infinity.
Bloom-Xu show that this is indeed the case \cite{bloom1} under the following hypothesis on the function $G$:

\medskip

\noindent (H4) For some $r_0 > 0$, we have
$$
\int_{r_0}^{\infty} r |G(r)| dr < +\infty
$$
and $G$ is bounded on $[r_0, \infty)$.

\medskip

Under hypothesis (H4), for any $\lambda \in \C - \{0\}$, there are unique solutions $\Phi_{\lambda}, \Phi_{-\lambda}$
of equation (\ref{evaleqn}) on $(0, \infty)$ which are asymptotic to exponentials near infinity \cite{bloom1},
$$
\Phi_{\pm \lambda}(r) = e^{(\pm i\lambda - \rho)r}(1 + o(1)) \ \hbox{ as } r \to +\infty
$$
The solutions $\Phi_{\lambda}, \Phi_{-\lambda}$ are linearly independent, so, since $\phi_{\lambda} = \phi_{-\lambda}$,
there exists a function $c$ on $\C - \{0\}$ such that
$$
\phi_{\lambda} = c(\lambda) \Phi_{\lambda} + c(-\lambda) \Phi_{-\lambda}
$$
for all $\lambda \in \C - \{0\}$. We will call this function the $c$-function of the hypergroup.
We remark that if the hypergroup $([0, \infty), *)$ is the one arising from
convolution of radial measures on a noncompact rank one symmetric space $G/K$, then this function agrees
with Harish-Chandra's $c$-function only on the half-plane $\{ \Im \lambda \leq 0 \}$ and not on all of $\C$.

\medskip

If we furthermore assume the hypothesis $|\alpha| \neq 1/2$, then Bloom-Xu show that the function $c$ is non-zero for
$\Im \lambda \leq 0, \lambda \neq 0$, and prove the following estimates:

\medskip

There exist constants $C, K > 0$ such that

\begin{align*}
\frac{1}{C} |\lambda| & \leq |c(\lambda)|^{-1} \leq C |\lambda| , \quad \quad  |\lambda|  \leq K \\
\frac{1}{C} |\lambda|^{\alpha + \frac{1}{2}} & \leq |c(\lambda)|^{-1} \leq C |\lambda|^{\alpha + \frac{1}{2}} , \quad |\lambda| \geq K \\
\end{align*}

Moreover they prove the following inversion formula: for any even function $f \in C^{\infty}_c(\R)$,
$$
f(r) = C_0 \int_{0}^{\infty} \hat{f}(\lambda) \phi_{\lambda}(r) |c(\lambda)|^{-2} d\lambda
$$
where $C_0 > 0$ is a constant.

\medskip

It follows that the Plancherel measure $\sigma$ of the hypergroup is supported on $[0, \infty)$, and absolutely
continuous with respect to Lebesgue measure, with density given by $C_0 |c(\lambda)|^{-2}$. Bloom-Xu also
show that the $c$-function is holomorphic on the half-plane $\{ \Im \lambda < 0 \}$.

\medskip

\subsection{The density function of a harmonic manifold}

\medskip

Let $X$ be a simply connected, $n$-dimensional harmonic manifold of purely exponential volume growth, and let $A$ be the density
function of $X$. We check that $A$ is a Chebli-Trimeche function, so that we obtain a commutative hypergroup $([0,\infty), *)$,
and that the conditions of Bloom-Xu are met so that the Plancherel measure is given by $C_0 |c(\lambda)|^{-2} d\lambda$ on
$[0, \infty)$.

\medskip

The function $A(r)$ equals, up to a constant factor, the volume of geodesic spheres $S(x, r)$, which is increasing in $r$ and
tends to infinity as $r$ tends to infinity, so condition (H1) is satisfied. As stated in section 2.2, the function $A'(r)/A(r)$
equals the mean curvature of geodesic spheres $S(x, r)$, which decreases
monotonically to a limit $h = 2\rho$ which is positive (and equals the mean curvature of horospheres), so condition (H2) is satisfied.

\medskip

Fixing a point $x \in X$, for $r > 0$, the density function $A(r)$ is given by the Jacobian of the map $\phi : v \mapsto \exp_x(rv)$ from
the unit tangent sphere $T^1_x X$ to the geodesic sphere $S(x, r)$. Let $T$ be the map $v \mapsto rv$ from the unit tangent sphere
$T^1_x X$ to the tangent sphere of radius $r$, $T^r_x X \subset T_x M$, then $\phi = \exp_x \circ T$, so the Jacobian of $\phi$ is
given by the product of the Jacobians of $T$ and $\exp_x$, hence
$$
A(r) = r^{n-1} B(r)
$$
where the function $B$ is given by
$$
B(r) = \det (D\exp_x)_{rv}
$$
where $v$ is any fixed vector in $T^1_x X$. Since $B$ is independent of the choice of $v$, in particular is the same for vectors $v$ and
$-v$, the function $B$ is even, and $C^{\infty}$ on $\R$ with $B(0) = 1$.
Thus condition (H3) holds for the function $A$, with $\alpha = (n - 2)/2$.

\medskip

The density function $A$ is thus a Chebli-Trimeche function, so we obtain a hypergroup structure on $[0, \infty)$, which we call
the {\it radial hypergroup} of the harmonic manifold $X$ (the reason for this terminology will become clear from the
the following sections).

\medskip

We proceed to check that condition (H4) is satisfied. For this we will need the following theorem of Nikolayevsy:

\medskip

\begin{theorem} \label{exppoly} \cite{nikolayevsky} The density function of a harmonic manifold is an {\it exponential polynomial},
i.e. a function of the form
$$
A(r) = \sum_{i = 1}^k (p_i(r) \cos(\beta_i r) + q_i(r) \sin(\beta_i r)) e^{\alpha_i r}
$$
where $p_i, q_i$ are polynomials and $\alpha_i, \beta_i \in \R$, $i = 1, \dots, k$.
\end{theorem}

\medskip

It will be convenient to rearrange terms and write the density function in the form
\begin{equation} \label{formofa}
A(r) = \sum_{i = 1}^l \sum_{j = 0}^{m_i} f_{ij}(r) r^j e^{\alpha_i r}
\end{equation}
where $\alpha_1 < \alpha_2 < \dots < \alpha_l$, and each $f_{ij}$ is a trigonometric polynomial,
i.e. a finite linear combination of functions of the form $\cos(\beta r)$ and $\sin(\beta r)$, $\beta \in \R$, with
$f_{im_i}$ not identically zero, for $i = 1, \dots, l$.
For an exponential polynomial written in this form, we will call the largest exponent $\alpha_l$ which appears in the
exponentials the {\it exponential degree} of the exponential polynomial.

\medskip

\begin{lemma} \label{densityform} With the density function as above, we have $\alpha_l = 2\rho, m_l = 0$ and $f_{l0} = C$
for some constant $C > 0$. Thus the density function is of the form
$$
A(r) = C e^{2\rho r} + P(r)
$$
where $P$ is an exponential polynomial of exponential degree $\delta < 2\rho$.
\end{lemma}

\medskip

\noindent{\bf Proof:} Recall that $X$ has purely exponential volume growth, i.e.
there exists a constant $C > 1$ such that
\begin{equation} \label{pureexp}
\frac{1}{C} \leq \frac{A(r)}{e^{2\rho r}} \leq C
\end{equation}
for all $r \geq 1$. If $\alpha_l < 2\rho$, then $A(r)/e^{2\rho r} \to 0$ as $r \to \infty$, contradicting (\ref{pureexp})
above, so we must have $\alpha_l \geq 2\rho$. On the other hand, if $\alpha_l > 2\rho$, then since
$f_{lm_l}$ is a trigonometric polynomial which is not identically zero, we can choose a sequence $r_m$ tending to infinity
such that $f_{lm_l}(r_m) \to \alpha \neq 0$. Then clearly $A(r_m)/e^{2\rho r_m} \to \infty$, again contradicting (\ref{pureexp}).
Hence $\alpha_l = 2\rho$.

\medskip

Using (\ref{formofa}) and $\alpha_l = 2\rho$, we have
$$
\frac{A'(r)}{A(r)} - 2\rho = \frac{f_{lm_l}'(r) + o(1)}{f_{lm_l}(r) + o(1)}
$$
as $r \to \infty$, thus
\begin{align*}
f_{lm_l}'(r) + o(1) & = (f_{lm_l}(r) + o(1))\left(\frac{A'(r)}{A(r)} - 2\rho\right) \\
                    & \to 0 \\
\end{align*}
as $r \to \infty$ since $f_{lm_l}$ is bounded and $A'(r)/A(r) - 2\rho \to 0$ as $r \to \infty$.
Thus $f_{lm_l}'$ is a trigonometric polynomial which tends to $0$ as $r \to \infty$, so it must be identically
zero, hence $f_{lm_l} = C$ for some non-zero constant $C$.

\medskip

It follows that
$$
A(r) = C r^{m_l} e^{2\rho r} (1 + o(1))
$$
as $r \to \infty$. If $m_l \geq 1$ then $A(r)/e^{2\rho r} \to \infty$ as $r \to \infty$, so we must have $m_l = 0$.
$\diamond$

\medskip

\begin{lemma} \label{h4holds} Condition (H4) holds for the density function $A$, i.e.
$$
\int_{r_0}^{\infty} r |G(r)| dr < +\infty
$$
and $G$ is bounded on $[r_0, \infty)$ for any $r_0 > 0$, where
$$
G(r) = \frac{1}{4} \left( \frac{A'(r)}{A(r)}\right)^2 + \frac{1}{2}\left(\frac{A'}{A}\right)'(r) - \rho^2
$$
\end{lemma}

\medskip

\noindent{\bf Proof:} By the previous lemma, $A(r) = C e^{2\rho r} + P(r)$, where $P$ is an
exponential polynomial of exponential degree $\delta < 2\rho$. We then have
\begin{align*}
\frac{A'(r)}{A(r)} - 2\rho & = \frac{P'(r) - 2\rho P(r)}{C e^{2\rho r} + P(r)} \\
                           & = \frac{Q(r)}{C e^{2\rho r} + P(r)} \\
\end{align*}
where $Q$ is an exponential polynomial of exponential degree less than or equal to $\delta$. Putting $\alpha = (2\rho - \delta)/2$,
it follows that $A'(r)/A(r) - 2\rho = O(e^{-\alpha r})$ as $r \to \infty$.
Differentiating, we obtain
\begin{align*}
\left(\frac{A'}{A}\right)'(r) & = \frac{(C e^{2\rho r} + P(r))Q'(r) - Q(r)(2\rho C e^{2\rho r} + P'(r))}{(C e^{2\rho r} + P(r))^2} \\
                 & = \frac{Q_1(r)}{(C e^{2\rho r} + P(r))^2} \\
\end{align*}
where $Q_1$ is an exponential polynomial of exponential degree less than or equal to $(2\rho + \delta)$. Since the denominator
of the above expression is of the form $k e^{4\rho r} + P_1(r)$ with $P_1$ an exponential polynomial of exponential degree strictly
less than $4\rho$, it follows that $(A'/A)'(r) = O(e^{-\alpha r})$ as $r \to \infty$.

\medskip

Now we can write the function $G$ as
$$
G(r) = \frac{1}{4}\left(\frac{A'(r)}{A(r)} - 2\rho\right)\left(\frac{A'(r)}{A(r)} + 2\rho\right) + \frac{1}{2}\left(\frac{A'}{A}\right)'(r)
$$
Since $(A'(r)/A(r) + 2\rho)$ is bounded, it follows from the previous paragraph that $G(r) = O(e^{-\alpha r})$ as $r \to \infty$.
This immediately implies that condition (H4) holds. $\diamond$

\medskip

In order to apply the result of Bloom-Xu on the Plancherel measure for the hypergroup, it remains to check
that $|\alpha| \neq 1/2$. Since $\alpha = (n - 2)/2$, this means $n \neq 3$. Now the
Lichnerowicz conjecture holds in dimensions $n \leq 5$ (\cite{lichnerowicz1}, \cite{walker1}, \cite{besse1}, \cite{nikolayevsky}), i.e. the only harmonic manifolds in such dimensions
are the rank one symmetric spaces $X = G/K$, for which as mentioned earlier the Jacobi analysis
applies, and the Plancherel measure of the hypergroup is well known to be given by $C_0 |{\bf c}(\lambda)|^{-2} d\lambda$ where
${\bf c}$ is Harish-Chandra's $c$-function. Thus in our case we may as well assume that $X$ has dimension $n \geq 6$, so that
$|\alpha| \neq 1/2$, and we may then apply the results of Bloom-Xu stated in the previous section.

\medskip

\subsection{The spherical Fourier transform}

\medskip

Let $\phi_{\lambda}$ denote as in section 4.1 the unique function on $[0, \infty)$ satisfying $L_R \phi_{\lambda} = -(\lambda^2 + \rho^2) \phi_{\lambda}$
and $\phi_{\lambda}(0) = 1$. For $x \in X$ let $d_x$ denote as before the distance function from the point $x$, $d_x(y) = d(x, y)$.
We define the following eigenfunction of $\Delta$ radial around $x$:
$$
\phi_{\lambda, x} := \phi_{\lambda} \circ d_x
$$
The uniqueness of $\phi_{\lambda}$ as an eigenfunction of $L_R$ with eigenvalue $-(\lambda^2 + \rho^2)$ and taking the value $1$ at $r = 0$
immediately implies the following lemma:

\medskip

\begin{lemma} \label{uniqueeigenfn} The function $\phi_{\lambda, x}$ is the unique eigenfunction $f$ of $\Delta$ on $X$ with
eigenvalue $-(\lambda^2 + \rho^2)$ which is radial around $x$ and satisfies $f(x) = 1$.
\end{lemma}

\medskip

Note that for $\lambda \in \R$, the functions $\phi_{\lambda, x}$ are bounded. Let $dvol$ denote the Riemannian volume measure on $X$.

\medskip

\begin{definition} Let $f \in L^1(X, dvol)$ be radial around the point $x \in X$. We define the spherical Fourier transform
of $f$ by
$$
\hat{f}(\lambda) := \int_X f(y) \phi_{\lambda, x}(y) dvol(y)
$$
for $\lambda \in \R$.
\end{definition}

\medskip

For $f$ a function on $X$ radial around the point $x$,
let $f = u \circ d_x$ where $u$ is a function on $[0, \infty)$, then evaluating the integral over $X$ in
geodesic polar coordinates gives
$$
\int_X |f(y)| dvol(y) = \int_{0}^{\infty} |u(r)| A(r) dr
$$
thus $f \in L^1(X)$ if and only if $u \in L^1([0, \infty), A(r) dr)$. In that case, again integrating in polar
coordinates gives
$$
\hat{f}(\lambda) = \int_{0}^{\infty} u(r) \phi_{\lambda}(r) A(r) dr = \hat{u}(\lambda)
$$
where $\hat{u}$ is the hypergroup Fourier transform of the function $u$. Moreover $f \in C^{\infty}_c(X)$ if and only if
$u$ extends to an even function on $\R$ such that $u \in C^{\infty}_c(\R)$. Applying the Fourier inversion formula of Bloom-Xu for the radial
hypergroup stated in section 4.1 to the function $u$ then leads immediately to the following inversion formula for radial functions:

\medskip

\begin{theorem} \label{radialinversion} Let $(X,g)$ be a simply connected harmonic manifold of purely exponential volume growth and $f \in C^{\infty}_c(X)$ be radial around the point $x \in X$. Then
$$
f(y) = C_0 \int_{0}^{\infty} \hat{f}(\lambda) \phi_{\lambda, x}(y) |c(\lambda)|^{-2} d\lambda
$$
for all $y \in X$. Here $c$ denotes the $c$-function of the radial hypergroup and $C_0 > 0$ is a constant.
Moreover, the $c$-function is holomorphic on the half-plane $\{ \Im \lambda < 0 \}$.
\end{theorem}

\medskip

\noindent{\bf Proof:} As shown in the previous section, all the hypotheses required to apply the inversion formula of Bloom-Xu are satisfied,
hence
$$
u(r) = C_0 \int_{0}^{\infty} \hat{u}(\lambda) \phi_{\lambda}(r) |c(\lambda)|^{-2} d\lambda
$$
Since $f = u \circ d_x$, this gives
\begin{align*}
f(y) & = u(d_x(y)) \\
     & = C_0 \int_{0}^{\infty} \hat{u}(\lambda) \phi_{\lambda}(d_x(y)) |c(\lambda)|^{-2} d\lambda \\
     & = C_0 \int_{0}^{\infty} \hat{f}(\lambda) \phi_{\lambda,x}(y) |c(\lambda)|^{-2} d\lambda \\
\end{align*}
For the holomorphicity of the function $c$ in $\{ \Im \lambda < 0 \}$ see the proof of Proposition 3.17 in \cite{bloom1}.
$\diamond$

\medskip

The Plancherel theorem for the radial hypergroup leads to the following:

\medskip

\begin{theorem} Let $(X,g)$ be a simply connected harmonic manifold of purely exponential volume growth. Let
$L^2_x(X, dvol)$ denote the closed subspace of $L^2(X)$ consisting of those functions
in $L^2(X)$ which are radial around the point $x$. For $f \in L^1(X, dvol) \cap L^2_x(X, dvol)$, we have
$$
\int_X |f(y)|^2 dvol(y) = C_0 \int_{0}^{\infty} |\hat{f}(\lambda)|^2 |c(\lambda)|^{-2} d\lambda
$$
The spherical Fourier transform $f \mapsto \hat{f}$ extends to an isometry from $L^2_x(X, dvol)$ onto
$L^2([0, \infty), C_0 |c(\lambda)|^{-2} d\lambda)$.
\end{theorem}

\medskip

\noindent{\bf Proof:} The map $u \mapsto f = u \circ d_x$ defines
an isometry of $L^2([0, \infty), A(r) dr)$ onto $L^2(X, dvol)_x$, which maps $L^1([0, \infty), A(r) dr) \cap L^2([0, \infty), A(r) dr)$
onto $L^1(X, dvol) \cap L^2_x(X, dvol)$. The statements of the theorem then follow from the Levitan-Plancherel theorem
for the radial hypergroup and from the fact that the Plancherel measure is supported on $[0, \infty)$,
given by $C_0 |c(\lambda)|^{-2} d\lambda$. $\diamond$

\medskip

\section{Fourier inversion and Plancherel theorem}

\medskip

As before, we assume in this section that $(X,g)$ denotes a simpy connected harmonic manifold of purely exponential volume growth unless stated otherwise. We proceed to the analysis of non-radial functions on $X$. Our definition of Fourier transform will depend on the choice of a basepoint $x \in X$.

\medskip

\begin{definition} Let $x \in X$. For $f \in C^{\infty}_c(X)$, the Fourier transform of $f$ based at the
point $x$ is the function
on $\C \times \partial X$ defined by
$$
\tilde{f}^x(\lambda, \xi) = \int_X f(y) e^{(-i\lambda - \rho)B_{\xi,x}(y)} dvol(y)
$$
for $\lambda \in \C, \xi \in \partial X$. Here as before $B_{\xi,x}$ denotes the Busemann function at
$\xi$ based at $x$ such that $B_{\xi,x}(x) = 0$.
\end{definition}

\medskip

Using the formula
$$
B_{\xi, x} = B_{\xi,o} - B_{\xi, o}(x)
$$
for points $o, x \in X$, we obtain the following relation between the Fourier transforms
based at two different basepoints $o, x \in X$:

\medskip

\begin{equation} \label{fourierbasept}
\tilde{f}^x(\lambda, \xi) = e^{(i\lambda + \rho)B_{\xi, o}(x)} \tilde{f}^o(\lambda, \xi)
\end{equation}

\medskip

The key to passing from the inversion formula for radial functions of section 4.3 to an inversion formula for non-radial
functions will be a formula expressing the radial eigenfunctions $\phi_{\lambda,x}$ as an integral with respect to
$\xi \in \partial X$ of the eigenfunctions $e^{(i\lambda - \rho)B_{\xi,x}}$ (Theorem \ref{phipoisson}). This will be the analogue of the well-known
formulae for rank one symmetric spaces $G/K$ and harmonic $NA$ groups expressing the radial eigenfunctions $\phi_{\lambda, x}$ as matrix
coefficients of representations of $G$ on $L^2(K/M)$ and $NA$ on $L^2(N)$ respectively. %for notational convenience, we will however denote the visibility measure on $\partial X$ by the same symbol $\lambda_x$.

\medskip

%In \cite{knieperpeyerimhoff2}, it is shown that the visibility measures $\lambda_x, x \in X$ are mutually absolutely
%continuous and their Radon-Nikodym derivatives are given by
%$$
%\frac{d\lambda_y}{d\lambda_x}(\xi) = e^{-2\rho B_{\xi, x}(y)}
%$$

We start with a basic relation between eigenfunctions of the Laplacian:

\begin{lemma} \label{radialisation} Let $x \in X$ and $\xi \in \partial X$. Then for all $\lambda \in \C$,
$$
\phi_{\lambda, x} = M_x (e^{(i\lambda - \rho)B_{\xi,x}})
$$
(where $M_x$ is the radialisation operator around the point $x$). In particular, $\phi_{\lambda, x}(y)$ is entire in
$\lambda$ for fixed $y \in X$, and is real and positive for $\lambda$ such that $(i\lambda - \rho)$ is real and positive.
\end{lemma}

\medskip

\noindent{\bf Proof:} Since the function $e^{(i\lambda - \rho)B_{\xi,x}}$ is an eigenfunction of the
Laplacian $\Delta$ with eigenvalue $-(\lambda^2 + \rho^2)$ and the operator $M_x$ commutes
with $\Delta$, the function $f = M_x (e^{(i\lambda - \rho)B_{\xi,x}})$ is also an eigenfunction of
$\Delta$ for the eigenvalue $-(\lambda^2 + \rho^2)$. Since $f$ is radial around $x$ and $f(x) = 1$, it follows
from Lemma \ref{uniqueeigenfn} that $f = \phi_{\lambda, x}$. $\diamond$

\medskip

The next proposition provides a connection between the Fourier transform and the spherical Fourier transform for radial functions:

\medskip

\begin{prop} \label{coincide} Let $f \in C^{\infty}_c(X)$ be radial around the point $x \in X$. Then the Fourier transform of $f$
based at $x$ coincides with the spherical Fourier transform,
$$
\tilde{f}^x(\lambda, \xi) = \hat{f}(\lambda)
$$
for all $\lambda \in \C, \xi \in \partial X$.
\end{prop}

\medskip

\noindent{\bf Proof:} Let $f = u \circ d_x$ where $u \in C^{\infty}_c(\R)$. By Lemma \ref{radialisation} above,
$$
\phi_{\lambda}(r) = \phi_{-\lambda}(r) = \int_{S(x,r)} e^{(-i\lambda - \rho)B_{\xi,x}(y)} d\sigma^r(y)
$$
where $\sigma^r$ is normalized surface area measure on the geodesic sphere $S(x,r)$.
Evaluating the integral defining $\tilde{f}^x$ in geodesic polar coordinates centered at $x$ we have
\begin{align*}
\tilde{f}^x(\lambda, \xi) & = \int_{0}^{\infty} \int_{S(x,r)} f(y) e^{(-i\lambda - \rho)B_{\xi,x}(y)} d\sigma^r(y) A(r) dr \\
                          & = \int_{0}^{\infty} u(r) \phi_{\lambda}(r) A(r) dr \\
                          & = \hat{f}(\lambda) \\
\end{align*}
$\diamond$

\medskip

Now we need to define the {\it visibility measures}
on the boundary $\partial X$: Given a point $x \in X$, let $\theta_x$ be normalized canonical measure on the unit tangent sphere $T^1_x X$,
i.e. the unique probability measure on $T^1_x X$ invariant under the orthogonal group of the tangent space
$T_x X$. For $v \in T^1_x X$, let $\gamma_v : [0, \infty) \to X$ be the unique geodesic ray with initial velocity $v$.
Then we have a homeomorphism $pr_x : T^1_x X \to \partial X, v \mapsto \gamma_v(\infty)$. The
visibility measure on $\partial X$ (with respect to the basepoint $x$)
is defined to be the push-forward $(pr_x)_* \theta_x$ of $\lambda_x$ under the map $pr_x$. 

\medskip

For $\lambda \in \C$ and $x \in X$, define the function $\tilde{\phi}_{\lambda, x}$ on $X$ by
$$
\tilde{\phi}_{\lambda,x}(y) = \int_{\partial X} e^{(i\lambda - \rho)B_{\xi, x}(y)} d\lambda_x(\xi)
$$
It follows from the above equation that $\tilde{\phi}_{\lambda, x}(y)$ is entire in $\lambda$ for fixed $y \in X$,
and is real and positive for $\lambda$ such that $(i\lambda - \rho)$ is real and positive. Moreover, by
Proposition \ref{helgasonkernel}, the function $\tilde{\phi}_{\lambda, x}$ is an eigenfunction of the Laplacian $\Delta$ with
eigenvalue $-(\lambda^2 + \rho^2)$, and $\tilde{\phi}_{\lambda, x}(x) = 1$.

Our next aim is to show that
$\tilde{\phi}_{\lambda,x}$ is radial around $x$ and, therefore, agrees with the function $\phi_{\lambda,x}$
introduced in Lemma \ref{uniqueeigenfn}. We start with a crucial property of non-compact harmonic manifolds without any further assumptions, derived from a result
of Szabo \cite{szabo} that the volume of the intersection of a metric ball $B(x,r_1)$ with a geodesic sphere
$S(y,r_2)$ depends only on the radii $r_1,r_2$ and the distance $d = d(x,y)$ of their centers. We will therefore denote this volume by $v(r_1,r_2,d)$.

\begin{prop}
  Let $(X,g)$ be a non-compact simply connected harmonic manifold. For $v \in T_x^1X$ and $r > 0$, let
  $b_v^r (y) = d(y,\gamma_v(r))-r$, and $\theta_x$ be the normalized canonical measure of $T^1_xX$. Then for every continuous function $\phi: \R \to \C$, the function
  $$ F(y) := \int_{T^1_xX} \phi(b_v^r(y)) d\theta_x(v) $$
  is radial around $x$.
\end{prop}

\medskip

\noindent{\bf Proof:}
  Let $\psi(s) = \phi(s-r)$. Then
  $$
  \phi(b_v^r(y)) = \phi(d(y,\gamma_v(r)) - r) = \psi(d(y,\gamma_v(r))
  $$
  and
  \begin{equation} \label{eq:F}
  F(y) = \int_{T	^1_xX} \phi(b_v^r(y)) d\theta_x(v) = \int_{T^1_xX} \psi(d(y,\gamma_v(r)) d\theta_x(v).
  \end{equation}
  Next, we consider the following expression:
  \begin{equation} \label{eq:fub}
  \int_{B(x,r)} \psi(d(y,z)) dvol(z) = \int_0^r A(t) \int_{T^1_xX} \psi(d(y,\gamma_v(t))) d\theta_x(v) dt.
  \end{equation}
  On the other hand, we have
  \begin{multline} \label{eq:re-interpret}
    \int_{B(x,r)} \psi(d(y,z)) dvol(z) = \int_0^\infty \int_{B(x,r) \cap S(y,t)} \psi(d(y,z)) d\sigma_{S_y(t)}(z) dt = \\
    \int_0^\infty \int_{B(x,r) \cap S(y,t)} \psi(t) d\sigma_{S_y(t)}(z) dt = \\
    \int_0^\infty vol_{S(y,t)}(B(x,r) \cap S(y,t)) \psi(t) dt = \int_0^\infty v(r,t,d(x,y)) \psi(t) dt.
   \end{multline}
   Now, we combine \eqref{eq:fub} and \eqref{eq:re-interpret} and differentiate with
   respect to $r$ and obtain
   $$ A(r) \int_{T^1_xX} \psi(d(y,\gamma_v(r)) d\theta_x(v) = \int_0^\infty \frac{\partial v}{\partial r}(r,t,d(x,y)) \psi(t) dt. $$
   In view of \eqref{eq:F}, this implies that
   $$ F(y) = \frac{1}{A(r)} \int_0^\infty \frac{\partial v}{\partial r}(r,t,d(x,y)) \psi(t) dt, $$
   which is obviously independent of the position of $y$ within the sphere $S(R,x)$ with $R = d(x,y)$.  This shows
   that the function $F$ is radial around $x$.
$\diamond$

\medskip

The analogous statement for Busemann functions is obtained via a limiting argument:

\medskip

\begin{cor} \label{cor:cruc}
  Let $(X,g)$ be a non-compact simply connected harmonic manifold and $\phi: \R \to \C$ be a continuous function. Then the function
  $$ F(y) := \int_{T^1_xX} \phi(b_v(y)) d\theta_x(v) $$
  is a radial function around $x$.
\end{cor}

\medskip

\noindent{\bf Proof:}
  Note that we have pointwise convergence $\phi(b_v^r(y)) \to \phi(b_v(y))$ for $r \to \infty$ and, since
  $$  |b_v^r(y)| \le d(x,y) \quad \text{for all $r \ge 0$,} $$
  we can apply Lebesgue's dominated convergence.
$\diamond$

\medskip

\begin{theorem} \label{phipoisson} Let $(X,g)$ be a non-compact simply connected harmonic manifold.
Let $\lambda \in \C$ and $x \in X$. Then
$$
\phi_{\lambda, x}(y) = \int_{T^1_x X} e^{(i\lambda - \rho)b_v(y)} d\theta_x(v)
$$
for all $y \in X$.
\end{theorem}

\medskip

\noindent{\bf Proof:}
  Both sides are eigenfunctions of the Laplacian $\Delta$ with eigenvalue $-(\lambda^2+\rho^2)$. Moreover, both sides
  assume the value $1$ as $y=x$. $\phi_{\lambda,x}$ is radial around $x$, by definition, and the right hand side
  is radial by Corollary \ref{cor:cruc} with $\phi(s) = e^{i\lambda-\rho)s}$. Therefore, both expressions agree by
  the uniqueness of radial solutions of $\Delta u = -(\lambda^2+\rho^2)u$, $u(x) = 1$.
$\diamond$

\medskip

We can now prove the Fourier inversion formula:

\medskip

\begin{theorem} \label{fourierinversion} Let $(X,g)$ be a simply connected harmonic manifold of purely
exponential volume growth. Fix a basepoint $o \in X$. Then for $f \in C^{\infty}_c(X)$ we have
$$
f(x) = C_0 \int_{0}^{\infty} \int_{\partial X} \tilde{f}^o(\lambda, \xi) e^{(i\lambda - \rho)B_{\xi, o}(x)} d\lambda_o(\xi) |c(\lambda)|^{-2} d\lambda
$$
for all $x \in X$ (where $C_0 > 0$ is a constant).
\end{theorem}

\medskip

\noindent{\bf Proof:} Given $f \in C^{\infty}_c(X)$ and $x \in X$, the function $M_x f$ is in $C^{\infty}_c(X)$, is radial
around the point $x$ and satisfies $(M_x f)(x) = f(x)$. By Theorem \ref{radialinversion} applied to the function $M_x f$ we have
\begin{align*} %\label{fx}
f(x) = (M_x f)(x) & = C_0 \int_{0}^{\infty} \widehat{M_x f}(\lambda) \phi_{\lambda, x}(x) |c(\lambda)|^{-2} d\lambda \\
                  & = C_0 \int_{0}^{\infty} \widehat{M_x f}(\lambda) |c(\lambda)|^{-2} d\lambda \\
\end{align*}
(since $\phi_{\lambda, x}(x) = 1$). Now using the formal self-adjointness of the operator $M_x$, Theorem \ref{phipoisson},
the fact that $\phi_{\lambda,x}$ is radial around $x$ and $\phi_{\lambda, x} = \phi_{-\lambda,x}$ we obtain
\begin{align*}
\widehat{M_x f}(\lambda) & = \int_X (M_x f)(y) \phi_{-\lambda, x}(y) dvol(y) \\
                     & = \int_X f(y) (M_x \phi_{-\lambda, x})(y) dvol(y) \\
                     & = \int_X f(y) \phi_{-\lambda, x}(y) dvol(y) \\
                     & = \int_X f(y) \left( \int_{T^1_x X} e^{(-i\lambda - \rho)b_v(y)} d\theta_x(v) \right) dvol(y) \\
                     & = \int_{T^1_x X} \left(\int_X f(y) \left( e^{(-i\lambda - \rho)b_v(y)} dvol(y)\right) d\theta_x(v) \right) \\
                     & = \int_{\partial X} \left(\int_X f(y) e^{(-i\lambda - \rho)B_{\xi, x}(y)} dvol(y)\right) d\lambda_x(\xi) \\
                     & = \int_{\partial X} \tilde{f}^x(\lambda, \xi) d\lambda_x(\xi) \\
\end{align*}

Using the relations \eqref{fourierbasept}, namely
$$
\tilde{f}^x(\lambda, \xi) = e^{(i\lambda + \rho)B_{\xi, o}(x)} \tilde{f}^o(\lambda, \xi)
$$
and \eqref{eq:RadNyk}, that is
$$
\frac{d\lambda_x}{d\lambda_o}(\xi) = e^{-2\rho B_{\xi, o}(x)},
$$
we get
\begin{align*}
\widehat{M_x f}(\lambda) & = \int_{\partial X} e^{(i\lambda + \rho)B_{\xi, o}(x)} \tilde{f}^o(\lambda, \xi) e^{-2\rho B_{\xi, o}(x)} d\lambda_o(\xi) \\
                     & = \int_{\partial X} \tilde{f}^o(\lambda, \xi) e^{(i\lambda - \rho)B_{\xi, o}(x)} d\lambda_o(\xi) \\
\end{align*}

Substituting this last expression for $\widehat{M_x f}(\lambda)$ in the equation 
$$f(x) = C_0 \int_{0}^{\infty} \widehat{M_x f}(\lambda) |c(\lambda)|^{-2} d\lambda$$
gives
$$
f(x) = C_0 \int_{0}^{\infty} \int_{\partial X} \tilde{f}^o(\lambda, \xi) e^{(i\lambda - \rho)B_{\xi, o}(x)} d\lambda_o(\xi) |c(\lambda)|^{-2} d\lambda
$$
as required.
$\diamond$

\medskip

The Fourier inversion formula leads immediately to a Plancherel theorem:

\medskip

\begin{theorem} \label{plancherel} Let $(X,g)$ be a simply connected harmonic manifold of purely exponential volume growth. Fix a basepoint $o \in X$. For $f, g \in C^{\infty}_c(X)$, we have
$$
\int_X f(x) \overline{g(x)} dvol(x) = C_0 \int_{0}^{\infty} \int_{\partial X} \tilde{f}^o(\lambda, \xi) \overline{\tilde{g}^o(\lambda, \xi)} d\lambda_o(\xi) |c(\lambda)|^{-2} d\lambda
$$
where $C_0$ is the constant appearing in the Fourier inversion formula.

\medskip

The Fourier transform $f \mapsto \tilde{f}^o$ extends to an isometry of $L^2(X, dvol)$ into
$L^2([0, \infty) \times \partial X, C_0 |c(\lambda)|^{-2} d\lambda d\lambda_o(\xi))$.
\end{theorem}

\medskip

\noindent{\bf Proof:} Applying the Fourier inversion formula to the function $g$ gives
\begin{align*}
\int_X f(x) & \overline{g(x)} dvol(x) \\ & =  C_0 \int_X f(x) \left(\int_{0}^{\infty} \int_{\partial X} \overline{\tilde{g}^o(\lambda, \xi)} e^{(-i\lambda - \rho)B_{\xi, o}(x)} d\lambda_o(\xi) |c(\lambda)|^{-2} d\lambda\right) dvol(x) \\
    & = C_0 \int_{0}^{\infty} \int_{\partial X} \left(\int_X f(x) e^{(-i\lambda - \rho)B_{\xi, o}(x)} dvol(x)\right) \overline{\tilde{g}^o(\lambda, \xi)} d\lambda_o(\xi) |c(\lambda)|^{-2} d\lambda \\
    & = C_0 \int_{0}^{\infty} \int_{\partial X} \tilde{f}^o(\lambda, \xi) \overline{\tilde{g}^o(\lambda, \xi)} d\lambda_o(\xi) |c(\lambda)|^{-2} d\lambda. 
\end{align*}
Taking $f = g$ gives that the Fourier transform preserves $L^2$ norms,
$$
||f||_2 = ||\tilde{f}^o||_2
$$
for all $f \in C^{\infty}_c(X)$. It follows from a standard argument that
the Fourier transform extends to an isometry of $L^2(X, dvol)$ into
$L^2([0, \infty) \times \partial X, C_0 |c(\lambda)|^{-2} d\lambda d\lambda_o(\xi))$.
$\diamond$

\medskip

\section{An integral formula for the $c$-function}

\medskip

In this section we prove the following identity which can be viewed as an analogue of a well-known integral formula for Harish-Chandra's {\bf c}-function (formula (18) in \cite{helgason1}, pg. 108):

\begin{theorem} \label{cfunction} Let $(X,g)$ be a simply connected harmonic manifold of purely exponential volume growth
and $c$ be the $c$-function of the radial hypergroup of $X$. Let $\Im \lambda < 0$. Then we have
$$
\lim_{r \to \infty} \frac{\phi_{\lambda}(r)}{e^{(i\lambda - \rho)r}} = c(\lambda) = \int_{\partial X} e^{-2(i\lambda - \rho)(\xi|\eta)_x} d\lambda_x(\eta).
$$
for any $x \in X, \xi \in \partial X$, where $(\xi|\eta)_x$ is the Gromov product given in Lemma \ref{def:gromov_product}.
\end{theorem}

For the proof of this identity we need some preparations.

\medskip

Recall that a geodesic metric space $(X,d)$ is called $\delta$-hyperbolic if geodesic triangles are $\delta$-thin,
that is each side is contained in the $\delta$-tubes of the other two sides. Moreover, the Gromov product $(y|z)_x$,
given by
$$ (y|z)_x = \frac{1}{2} (d(x,y)+d(x,z)-d(y,z)), $$
satisfies the following straightforward consequence of the triangle inequality: Let $\gamma$ be a geodesic joining
$y,z \in X$. Then for any point $w$ on this geodesic $\gamma$ we have
$$ (y,z)_x \le d(x,w). $$
This inequality entends to the boundary:
$$ (\xi|\eta)_x \le d(x,w), $$
for all points $w$ on any geodesic connecting $\xi, \eta \in \partial X$.

We use the Gromov product to define
balls in the boundary $\partial X$ with center $\xi \in \partial X$ and radius $r > 0$:
$$ B^{(x)}(\xi,r) := \{ \eta \in \partial X \mid e^{-(\xi|\eta)_x} < r \}. $$
Note that these ''balls'' do not come from a metric but from the Gromov product. We need the following geometric result.

\begin{lemma}
  Let $(X,d)$ be $\delta$-hyperbolic, $x \in X$ and $\gamma_{x,\xi}: [0,\infty) \to X$ be a geodesic ray with $\gamma_{x,\xi}(0) = x$ and  $\gamma_{x,\xi}(\infty) = \xi \in \partial X$. Then we have for all  $\epsilon \in (0,1)$, $y=\gamma_{x,\xi}(\log(1/\epsilon))$ and all $\eta \in B^{(x)}(\xi,\epsilon)$:
  \begin{equation} \label{eq:busestimate}
  \vert B_{\eta,y}(x) - d(x,y) \vert \le 6 \delta.
  \end{equation}
\end{lemma}

\medskip

\noindent {\bf Proof:} Let $\eta \in B^{(x)}(\xi,\epsilon)$ be fixed and $R = \log(1/\epsilon)$. Then $(\xi|\eta)_x \ge R$. Let
$\gamma_{\xi,\eta}: \R \to X$ be a geodesic connecting $\xi$ and $\eta$ and $\gamma_{x,\eta}: [0,\infty) \to X$ be a geodesic ray
connecting $x$ and $\eta$. Let $y_0 = \gamma_{x,\xi}(R - 2\delta)$. Then $y_0$ is not contained in the
$\delta$-tube around $\gamma_{\xi,\eta}(\R)$ since $d(x,y_0) = R - 2\delta$ and $d(x,\gamma_{\xi,\eta}(\R)) \ge (\xi|\eta)_x \ge R$.
Since triangles are $\delta$-thin, $y_0$ is contained in the $\delta$-tube around $\gamma_{x,\eta}(0,\infty)$. Let
$z_0 \in \gamma_{x,\eta}(0,\infty)$ with $d(y_0,z_0) \le \delta$ and, therefore, $d(y,z_0) \le 3 \delta$. This implies
for $z = \gamma_{x,\eta}(t)$ and $t > 0$ large:
\begin{multline*}
| d(x,z) - d(y,z) - d(x,y) | \le \\
%| d(x,z) - d(z_0,z) + (d(z_0,z) - d(y,z)) - d(x,z_0) + (d(x,z_0)-d(x,y)) | \le \\
| d(x,z) - d(z_0,z) - d(x,z_0) | + | d(z_0,z) - d(y,z) | + | d(x,z_0)-d(x,y) | \le 6 \delta
\end{multline*}
since $x,z_0,z$ lie on the geodesic $\gamma_{x,\eta}$ and, therefore, $d(x,z) - d(z_0,z) - d(x,z_0) = 0$ and $| d(z_0,z) - d(y,z) | , | d(x,z_0)-d(x,y) | \le d(y,z_0) \le 3 \delta$. The result follows then by taking the limit $t \to \infty$.
$\diamond$

\medskip

This result has the following consequence:

\begin{lemma} \label{lem:visualmeasurebound}
  Let $(X,g)$ be a non-compact simply connected $\delta$-hyperbolic harmonic manifold with horospheres of mean curvature $h > 0$. Then
  we have for all $x \in X$, $\xi \in \partial X$
  and $\epsilon \in (0,1)$:
  $$ \lambda_x(B^{(x)}(\xi,\epsilon)) \le e^{6\delta h} \epsilon^h. $$
\end{lemma}

%\begin{rmk}
%  This Lemma includes the case of a simply connected harmonic manifold of purely exponential volume growth since the conditions
%   $\Vert R \Vert \le R_0$ and $\Vert \nabla R \Vert \le R_0'$ are satisfied (\cite[Propositions 6.57 and 6.68]{besse1}) and purely exponential volume growth is equivalent to Gromov hyperbolicity (\cite{knieperpeyerimhoff3}).
%\end{rmk}

\noindent {\bf Proof:}  Recall that Gromov hyperbolicity and purely exponential volume growth are equivalent in the setting of non-compact simply connected harmonic manifolds (\cite{knieper1}). We use \cite[Theorem 1.4]{knieperpeyerimhoff2} (see also \eqref{eq:RadNyk}) about the Radon-Nykodym derivative and \eqref{eq:busestimate} to obtain for $y = \gamma_{x,\xi}(\log(1/\epsilon))$ with $\gamma_{x,\xi}$ a geodesic ray connecting $x$ and $\xi$:
\begin{multline*}
\lambda_x(B^{(x)}(\xi,\epsilon)) = \int_{B^{(x)}(\xi,\epsilon)} d\lambda_x(\eta) = \int_{B^{(x)}(\xi,\epsilon)}
\frac{d\lambda_x}{d\lambda_y} d\lambda_y(\eta) = \\ \int_{B^{(x)}(\xi,\epsilon)} e^{-h B_{\eta,y}(x)} d\lambda_y(\eta)
= \int_{B^{(x)}(\xi,\epsilon)} e^{-h B_{\eta,y}(x)} d\lambda_y(\eta) = \\ \int_{B^{(x)}(\xi,\epsilon)} e^{-h d(x,y)} e^{-h (B_{\eta,y}(x) - d(x,y))} d\lambda_y(\eta) \leq \epsilon^h \int_{B^{(x)}(\xi,\epsilon)} e^{6 \delta h} d\lambda_y(\eta) = e^{6 \delta h} \epsilon^h.
\end{multline*}
$\diamond$

 \medskip

With these results we can now present the proof of Theorem \ref{cfunction}:

\medskip

\noindent{\bf Proof:} For $\Im \lambda < 0$, using $\phi_{\lambda} = c(\lambda) \Phi_{\lambda} + c(-\lambda) \Phi_{-\lambda}$ and
$$ \Phi_{\pm \lambda}(r) = e^{(\pm i \lambda - \rho)r}(1 + o(1)) \quad \text{as $r \to \infty$}, $$
we have
\begin{align} \label{part1}
\frac{\phi_{\lambda}(r)}{e^{(i\lambda - \rho)r}} & = c(\lambda)(1 + o(1)) + c(-\lambda)e^{-2i\lambda r}(1 + o(1)) \\
                                                 & \to c(\lambda) \nonumber
\end{align}
as $r \to \infty$. This proves the first equation in the theorem.

\medskip

For the second equation in the theorem, we first consider the case $\lambda = it$ where $t \leq -\rho$, so that $\mu := i\lambda - \rho \geq 0$. Fix $x \in X$ and $\xi \in \partial X$.
For $\eta \in \partial X$, let $\gamma_{x,\eta}: [0,\infty) \to X$ be the geodesic ray satisfying $\gamma_{x,\eta}(0) = x$ and
$\gamma_{x,\eta}(\infty) = \eta$. The normalized surface area measure on the geodesic sphere $S(x, r)$ is given by the
push-forward of $\lambda_x$ under the map $\eta \mapsto \gamma_{x,\eta}(r)$, so by Lemma \ref{radialisation}
$$
\frac{\phi_{\lambda}(r)}{e^{(i\lambda - \rho)r}} = \int_{\partial X} e^{(i\lambda - \rho)(B_{\xi, x}(\gamma_{x,\eta}(r)) - r)} d\lambda_x(\eta)
$$

We will apply the dominated convergence theorem to evaluate the limit of the above integral as $r \to \infty$. First note that
by Lemma \ref{buseform}, for any $\eta$ not equal to $\xi$,
$$
B_{\xi, x}(\gamma_{x,\eta}(r)) - r \to -2(\xi|\eta)_x
$$
as $r \to \infty$, so the integrand converges a.e. as $r \to \infty$,
$$
e^{(i\lambda - \rho)(B_{\xi, x}(y(\eta, r)) - r)} \to e^{-2(i\lambda - \rho)(\xi|\eta)_x}.
$$
Now, using $|B_{\xi, x}(\gamma_{x,\eta}(r))| \le d(x,\gamma_{x,\eta}(r)) = r$ and $\mu \geq 0$ we have
$$
e^{\mu (B_{\xi, x}(\gamma_{x,\eta}(r)) - r)} \le 1.
$$
So dominated convergence applies and we conclude that
$$
\frac{\phi_{\lambda}(r)}{e^{(i\lambda - \rho)r}} \to \int_{\partial X} e^{-2(i\lambda - \rho)(\xi|\eta)_x} d\lambda_x(\eta)
$$
as $r \to \infty$. This shows the equation
$$ c(\lambda) = \int_{\partial X} e^{-2(i\lambda - \rho)(\xi|\eta)_x} d\lambda_x(\eta) $$
for $\lambda = it, t \leq -\rho$. Since $c(\lambda)$ is holomorphic for $\Im \lambda < 0$, we need to show
that the right hand side is also holomorphic for $\Im \lambda <0$. Then both expressions must be equal
for $\Im \lambda <0$, finishing the proof of the theorem.

Since $e^{-2(i\lambda - \rho)(\xi|\eta)_x}$ is holomorphic for all $\lambda \in \C$, we need to show that
$$ \int_{\partial X} \vert e^{-2(i\lambda - \rho)(\xi|\eta)_x} \vert d\lambda_x(\eta) < \infty $$
for $\Im \lambda < 0$. Then this expression is holomorphic for $\Im \lambda < 0$ by Morera's Theorem.
Let $\lambda = \sigma - i\tau$ with $\sigma \in \R$ and $\tau > 0$. Then we have
\begin{eqnarray*}
\int_{\partial X} \vert e^{-2(i\lambda - \rho)(\xi|\eta)_x} \vert d\lambda_x(\eta) &=& \int_{\partial X}
e^{-2 (\tau-\rho) (\xi|\eta)_x} d\lambda_x(\eta) \\
&=& \int_0^\infty \lambda_x( \{ \eta \in \partial X \mid e^{-2 (\tau-\rho) (\xi|\eta)_x} > t \} ) dt.
\end{eqnarray*}
If $\tau \ge \rho$ then the set $\{ \eta \in \partial X \mid e^{-2 (\tau-\rho) (\xi|\eta)_x} > t \}$ is empty
for $t > 1$, and so the last integral reduces to an integral over $[0,1]$, which is bounded above
by one since $\lambda_x$ is a probability measure.

Since $X$ is of purely exponential volume growth, it is a $\delta$-hyperbolic space for some $\delta > 0$ (\cite{knieper1}).
For $0 < \tau < \rho$ using Lemma
\ref{lem:visualmeasurebound} and the fact that $\lambda_x$ is a probability
measure we obtain with $h = 2\rho$
\begin{eqnarray*}
\int_0^\infty \lambda_x( \{ \eta \mid e^{-2 (\tau-\rho) (\xi|\eta)_x} > t \} ) dt &\le& 1 + \int_1^\infty \lambda_x( B^{(x)}(\xi,
( 1/t )^{1/(2(\rho-\tau))})) dt \\
&\le& 1 + e^{6\delta h} \int_1^\infty \left( \frac{1}{t} \right)^{\frac{2\rho}{2(\rho-\tau)}} dt \\
&<& \infty.
\end{eqnarray*}
$\diamond$

\section{The convolution algebra of radial functions}

In this section, we assume $(X,g)$ to be a non-compact simply connected harmonic manifold without
any further assumption unless stated otherwise. Fix a basepoint $o \in X$. We define a notion of convolution with radial functions as follows:

\medskip

For a function $f$ radial around the point $o$, let $f = u \circ d_o$, where $u$ is a function on $\R$. For $x \in X$,
the {\it x-translate} of $f$ is defined to be the function
$$
\tau_x f = u \circ d_x
$$
Note that if $f \in L^1(X, dvol)$, then evaluating integrals in geodesic polar coordinates centered at $o$ and $x$ gives
$$
||f||_1 = \int_{0}^{\infty} |u(r)| A(r) dr = ||\tau_x f||_1
$$

\begin{definition} \label{convolution} For $f$ an $L^1$ function on $X$ and $g$ an $L^1$ function on $X$ which is radial
around the point $o$, the convolution of $f$ and $g$ is the function on $X$ defined by
$$
(f * g)(x) = \int_X f(y) (\tau_x g)(y) dvol(y)
$$
\end{definition}

\medskip

Note that, if $g = u \circ d_o$, then
\begin{align*}
||f * g ||_1 & \leq \int_X \int_X |f(y)| |(\tau_x g)(y)| dvol(y) dvol(x) \\
             & = \int_X |f(y)| \left(\int_X |u(d(x, y))| dvol(x)\right) dvol(y) \\
             & = \int_X |f(y)| \left(\int_{0}^{\infty} |u(r)| A(r) dr\right) dvol(y) \\
             & = ||f||_1 ||g||_1 \\
             & < +\infty \\
\end{align*}
so that the integral defining $(f * g)(x)$ exists for a.e. x, and $f * g \in L^1(X, dvol)$.

\medskip

\begin{theorem} \label{radialconvolution} Let $(X,g)$ be a non-compact simply connected harmonic manifold. Let $L^1_o(X, dvol)$ denote the closed subspace of
$L^1(X, dvol)$ consisting of those
$L^1$ functions which are radial around the point $o$. Then for $f, g \in L^1_o(X, dvol)$ we have $f * g \in L^1_o(X, dvol)$,
and $L^1_o(X, dvol)$ forms a commutative Banach algebra under convolution.
\end{theorem}

\medskip

\noindent{\bf Proof:} We first consider functions $f, g \in C^{\infty}_c(X)$ which are radial around $o$.
It was shown in \cite[Lemma 2.8]{PS15} that $f * g$ is again radial around $o$ and it follows from \cite[Remark 1, p.127]{PS15} that $f*g=g*f$.

%For $x \in X$, the function $M_x \phi_{\lambda, o}$ is an eigenfunction of $\Delta$ with eigenvalue $-(\lambda^2 + \rho^2)$
%which is radial around the point
%$x$ and takes the value $\phi_{\lambda, o}(x)$ at the point $x$, so it follows from Lemma \ref{uniqueeigenfn} that

%$$
%M_x \phi_{\lambda, o} = \phi_{\lambda, o}(x) \phi_{\lambda, x}
%$$

%The above equation together with the spherical Fourier inversion formula for the function $f$ gives, for $y \in Y$,
%\begin{align*}
%(M_x f)(y) & = C_0 \int_{0}^{\infty} (M_x \phi_{\lambda, o})(y) \hat{f}(\lambda) |c(\lambda)|^{-2} d\lambda \\
%           & = C_0 \int_{0}^{\infty} \phi_{\lambda, o}(x) \phi_{\lambda, x}(y) \hat{f}(\lambda) |c(\lambda)|^{-2} d\lambda \\
%\end{align*}

%Now using the formal self-adjointness of $M_x$ we have
%\begin{align*}
%(f * g)(x) & = \int_X f(y) (\tau_x g)(y) dvol(y) \\
%           & = \int_X f(y) (M_x \tau_x g)(y) dvol(y) \\
%           & = \int_X (M_x f)(y) (\tau_x g)(y) dvol(y) \\
%           & = C_0 \int_{0}^{\infty} \phi_{\lambda, o}(x) \hat{f}(\lambda) \left(\int_X  \phi_{\lambda, x}(y) (\tau_x g)(y) dvol(y)\right) |c(\lambda)|^{-2} d\lambda \\
%           & = C_0 \int_{0}^{\infty} \phi_{\lambda, o}(x) \hat{f}(\lambda) \hat{g}(\lambda) |c(\lambda)|^{-2} d\lambda \\
%\end{align*}

%It follows from the above equation that $f * g$ is radial around the point $o$ since all the functions $\phi_{\lambda, o}$ are radial
%around the point $o$, and it also follows that $f * g = g * f$.

%\medskip

Now the inequality $||f * g||_1 \leq ||f||_1 ||g||_1$ implies, by the density of smooth, compactly supported radial functions in the space
$L^1_o(X, dvol)$, that for $f, g \in L^1_o(X, dvol)$ we have $f * g = g * f \in L^1_o(X, dvol)$, so $L^1_o(X, dvol)$ forms a
commutative Banach algebra under convolution. $\diamond$

\medskip

Now we derive a basic identity about the Fourier transform of a convolution. We assume here additionally that $(X,g)$ is of purely exponential volume growth to guarantee the existence of the Fourier transform. Note if $f, g \in C^{\infty}_c(X)$ with $g = u \circ d_o$ radial around $o$, then $f * g$ is compactly supported. For the
Fourier transform of $f * g$ based at $o$, using the identity $B_{\xi,o}(x) = B_{\xi,o}(y) + B_{\xi, y}(x)$ we have
\begin{align*}
\widetilde{f * g}^o(\lambda, \xi) & = \int_X \left(\int_X f(y)u(d(x,y)) dvol(y)\right) e^{(-i\lambda - \rho)B_{\xi, o}(x)} dvol(x) \\
                                  & = \int_X f(y) e^{(-i\lambda - \rho)B_{\xi, o}(y)} \left(\int_X u(d(x,y)) e^{(-i\lambda - \rho)B_{\xi, y}(x)} dvol(x)\right) dvol(y) \\
                                  & = \int_X f(y) e^{(-i\lambda - \rho)B_{\xi, o}(y)} \widetilde{u \circ d_y}^y(\lambda, \xi) dvol(y) \\
                                  & = \tilde{f}^o(\lambda, \xi) \hat{u}(\lambda) \\
                                  & = \tilde{f}^o(\lambda, \xi) \hat{g}(\lambda) \\
\end{align*}
where we have used the fact that for the function $u \circ d_y$ which is radial around $y$ we have
$$
\widetilde{u \circ d_y}^y(\lambda, \xi) = \hat{u}(\lambda) = \hat{g}(\lambda)
$$
where $\hat{u}$ is the hypergroup Fourier transform of $u$ and $\hat{g}$ is the spherical Fourier transform of the function
$g$ which is radial around $o$.

\medskip

Finally, we remark that the radial hypergroup of a harmonic manifold $(X,g)$ of purely exponential volume growth can be realized as the convolution algebra
of finite radial measures on the manifold: convolution with radial measures can be defined, and the convolution of two radial
measures is again a radial measure. This can be proved by approximating finite radial measures by $L^1$ radial functions and
applying the Theorem \ref{radialconvolution}. The convolution algebra $L^1_o(X, dvol)$ is then identified with a subalgebra of the hypergroup algebra
of finite radial measures under convolution.

\medskip

\medskip

\section{The Kunze-Stein phenomenon}

\medskip

In this section we assume that $(X,g)$ is a simply connected harmonic manifold of purely exponential volume growth and we prove a version of the Kunze-Stein phenomenon: for $1 \leq p < 2$, convolution with a radial $L^p$-function
defines a bounded operator on $L^2(X)$.

\medskip

\begin{lemma} \label{pstrip} Let $x \in X$, let $q > 2$, and let $\gamma_q = 1 - \frac{2}{q}$. Then for any $t \in (-\gamma_q \rho, \gamma_q \rho)$,
for any $\lambda \in \mathbb{C}$ with $\Im \lambda = t$ we have
$$
||\phi_{\lambda, x}||_q \leq ||\phi_{it, x}||_q < +\infty
$$
\end{lemma}

\medskip

\noindent{\bf Proof:} Given $t \in (-\gamma_q \rho, \gamma_q \rho)$, by Theorem \ref{phipoisson}, for $\lambda$ with $\Im \lambda = t$, we have
for any $y \in X$,
\begin{align*}
|\phi_{\lambda, x}(y)| & = \left| \int_{\partial X} e^{(i\lambda - \rho)B_{\xi, x}(y)} d\lambda_x(\xi) \right| \\
                       & \leq \int_{\partial X} e^{(-t - \rho)B_{\xi, x}(y)} d\lambda_x(\xi) \\
                       & = \phi_{it, x}(y) \\
\end{align*}
hence
$$
||\phi_{\lambda, x}||_q \leq ||\phi_{it, x}||_q
$$
If $t \neq 0$, then since $\phi_{it, x} = \phi_{-it, x}$, we may as well assume that $t > 0$, in which case we have,
letting $r = d(x,y)$,
\begin{align*}
\phi_{it, x}(y) & = c(it)\Phi_{it}(r) + c(-it)\Phi_{-it}(r) \\
                & = c(it) e^{(-t-\rho)r}(1 + o(1)) + c(-it) e^{(t - \rho)r}(1 + o(1)) \\
                & = c(-it) e^{(t - \rho)r}(1 + o(1)) \\
\end{align*}
as $r \to \infty$, so $|\phi_{it, x}(y)| \leq C e^{(t - \rho)r}$ for $r \geq M$ for some constants $C, M > 0$. We may also assume
$A(r) \leq C e^{2\rho r}$ for $r \geq M$. Then, evaluating integrals in geodesic polar coordinates centered at $x$, we have
\begin{align*}
\int_{d(x, y) \geq M} |\phi_{it,x}(y)|^q dvol(y) & \leq \int_{M}^{\infty} (C e^{(t - \rho)r})^q (C e^{2\rho r}) dr \\
                                                 & < +\infty \\
\end{align*}
since $(t - \rho)q + 2\rho < 0$ for $0 < t < \gamma_q \rho$, thus $||\phi_{it, x}||_q < +\infty$.

\medskip

For $t = 0$, applying H{\"o}lder's inequality we have, for any $\epsilon > 0$,
\begin{align*}
\phi_{0, x}(y) & = \int_{\partial X} e^{-\rho B_{\xi, x}(y)} d\lambda_x(\xi) \\
               & = \left( \int_{\partial X} e^{-(1 + \epsilon) \rho B_{\xi, x}(y)} d\lambda_x(\xi) \right)^{1/(1 + \epsilon)} \\
               & = \phi_{i\epsilon, x}(y)^{1/(1 + \epsilon)} \\
\end{align*}
from which it follows that by choosing $\epsilon$ small enough so that $q/(1 + \epsilon) > 2$ we have $||\phi_{0, x}||_q < +\infty$.
$\diamond$

\medskip

We remark that while the spherical Fourier transform was originally defined for radial $L^1$ functions, after fixing
a basepoint $x \in X$ it can also be defined for general $L^1$ functions by the same formula
$$
\hat{g}(\lambda) := \int_X g(y) \phi_{\lambda, x}(y) dvol(y) \ , \ \lambda \in \R
$$

We then have the following Lemma:

\medskip

\begin{lemma} \label{lpfourier} Let $x \in X$, let $1 \leq p < 2$ and let $g$ be an $L^p$-function on $X$. Let $q > 2$ be such that $\frac{1}{p} + \frac{1}{q} = 1$. Then the spherical Fourier transform $\hat{g}$ of $g$ extends to a holomorphic function of $\lambda$ on the
strip $S_q := \{ |\Im \lambda| < \gamma_q \rho \}$, and is bounded on any closed sub-strip $\{ |\Im \lambda| \leq t \}$ for $0 < t < \gamma_q \rho$.
In particular $\hat{g}$ on $\R$ satisfies a bound
$$
||\hat{g}||_{\infty} \leq C_p ||g||_p
$$
for a constant $C_p > 0$.
\end{lemma}

\medskip

\noindent{\bf Proof:} Given $0 < t < \gamma_q \rho$, for any $\lambda \in \C$ with $|\Im \lambda| \leq t$, by the previous Lemma
$||\phi_{\lambda, x}||_q \leq C$ for some constant $C$ only depending on $q$ and $t$, so it follows from Holder's inequality that the
function
$$
\hat{g}(\lambda) := \int_X g(y) \phi_{\lambda, x}(y) dvol(y)
$$
is well-defined and bounded for $|\Im \lambda| \leq t$ by a constant $C_{q,t}$ times $||g||_p$.
The holomorphicity of the function $\hat{g}$ follows from Morera's theorem,
using the holomorphic dependence of $\phi_{\lambda,x}$ on $\lambda$. $\diamond$

\medskip

We can now prove the following version of the Kunze-Stein phenomenon:

\medskip

\begin{theorem} \label{kunzestein} Let $(X,g)$ be a simply connected harmonic manifold of purely exponential volume growth. Let $x \in X$ and let $1 \leq p < 2$. Let $g \in C^{\infty}_c(X)$ be
radial around the point $x$. Then for any $f \in C^{\infty}_c(X)$ we have
$$
||f * g||_2 \leq C_p ||g||_p ||f||_2
$$
for some constant $C_p > 0$.
It follows that for any $g \in L^p(X)$ radial around $x$, the map $f \in C^{\infty}_c(X) \mapsto f * g$
extends to a bounded linear operator on $L^2(X)$ with operator norm at most $C_p ||g||_p$.
\end{theorem}

\medskip

\noindent{\bf Proof:} Recall that for $f, g \in C^{\infty}_c(X)$ with $g$ radial around $x$,
the Fourier transform of a convolution satisfies
$$
\widetilde{f * g}^x(\lambda, \xi) = \tilde{f}^x(\lambda, \xi) \hat{g}(\lambda)
$$
for $\lambda \in \R, \xi \in \partial X$. Applying the Plancherel theorem and Lemma \ref{lpfourier} above, we have
\begin{align*}
||f * g||_2 & = ||\widetilde{f * g}^x||_2 \\
            & = ||\tilde{f}^x \hat{g}||_2 \\
            & \leq ||\hat{g}||_{\infty} ||\tilde{f}^x||_2 \\
            & \leq C_p ||g||_p ||f||_2 \\
\end{align*}
The above inequality, valid for $C^{\infty}_c$-functions, implies by a standard density argument that for any
$L^p$ radial function $g$, the map $f \in C^{\infty}_c(X) \mapsto f * g$ extends to a bounded linear operator on
$L^2(X)$ with norm at most $C_p ||g||_p$. $\diamond$

\bibliography{moeb}
\bibliographystyle{alpha}

\end{document}